\documentclass[11pt]{article}
\usepackage{times}
\usepackage{amsfonts}
\usepackage{amsmath}
\usepackage{latexsym}
\usepackage{epsfig}
\usepackage{epsf}
\usepackage{epsfig, amsmath, amsthm}
\newtheorem{df}{Definition}[section]

\newtheorem{theorem}[df]{Theorem}
\newtheorem{lemma}[df]{Lemma}

\newtheorem{claim}[df]{Claim}
\newtheorem{corollary}[df]{Corollary}

\newtheorem{fact}[df]{Fact}

\begin{document}

\title{ On Computing  the Distinguishing Numbers of \\ Planar Graphs and Beyond:  a Counting Approach\footnote{A preliminary version of this paper~\cite{AD} appeared in the Proceedings of the Nordic Combinatorial Conference in 2004.}}
\author{V. Arvind \\ The Institute of Mathematical Sciences \\ Chennai 600 113, India \\ {\tt arvind@imsc.res.in} \and Christine T. Cheng \\ Department of Electrical Engineering and Computer Science \\ University of Wisconsin-Milwaukee \\ Milwaukee, WI 53211 \\ {\tt ccheng@uwm.edu} \and Nikhil R. Devanur \\ College of Computing \\ Georgia Institute of Technology \\ Atlanta, GA 30332 \\ {\tt nikhil@cc.gatech.edu}}
\date{\today}
\maketitle

\abstract{A vertex $k$-labeling of graph $G$ is  {\it distinguishing} if the only automorphism that preserves the labels of $G$ is the identity map.  The {\it distinguishing number of $G$}, $D(G)$, is the smallest integer $k$ for which $G$ has a distinguishing $k$-labeling.  In this paper, we  apply the principle of inclusion-exclusion and develop recursive formulas to count the number of inequivalent distinguishing $k$-labelings of a graph.  Along the way, we prove that the distinguishing number of a planar graph can be computed in time polynomial in the size of the graph.}

\section{Introduction}


 A vertex $k$-labeling of graph $G$ is a mapping $\phi: V(G) \rightarrow \{1, 2, \hdots, k\}$.  It is said to be {\it distinguishing} if the only automorphism that preserves the labels of $G$ is the identity map.  The {\it distinguishing number} of $G$, $D(G)$, is the minimum number of labels needed so that $G$ has a distinguishing labeling. The notion of distinguishing numbers for graphs was first introduced and developed by Albertson and Collins~\cite{AC}.  Their focus was on determining the relationships between a graph's automorphism group and its distinguishing number.  Their work has since been extended in many directions by researchers for graphs and groups (e.g., \cite{Albertson,Cowen,Chan2, Chan1, Chan3, CT, Tucker2, IK, klavzar, Tucker1, Julianna}). 

\begin{figure}
\label{compare-fig}
\centering
\epsfig{file=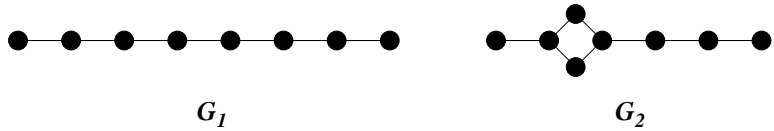, width=4in} 
\caption{An example of two graphs with the same number of vertices, isomorphic automorphism groups and identical distinguishing numbers but different number of inequivalent $k$-labelings: $D(G_1,k) =  k^4(k^4-1)/2$ but $D(G_2,k) = k^7(k-1)/2$.}
\end{figure}

Let $(G,\phi)$ denote the labeled version of $G$ under the labeling $\phi$. 
Given two distinguishing $k$-labelings $\phi$ and $\phi'$ of $G$, we say that $\phi$ and $\phi'$ are {\it equivalent} if there is some automorphism of $G$ that maps $(G,\phi)$ to $(G,\phi')$. 
We are interested in computing $D(G,k)$ -- the number of inequivalent $k$-distinguishing labelings of $G$ -- which was first considered by Arvind and Devanur~\cite{AD} and Cheng~\cite{Cheng} to determine the distinguishing numbers of trees.  Our motivation for studying this parameter are as follows.  First, $D(G) = \min\{k : D(G,k) > 0 \}$ so if we can compute  $D(G,k)$ efficiently then we can also determine $D(G)$ efficiently.  The usual way of proving that $D(G) = k^*$ is to present a distinguishing $k^*$-labeling of $G$ and then  argue that $G$ has no distinguishing labelings that uses $k < k^*$ labels.  Counting the number of  inequivalent distinguishing $k$-labelings of $G$ provides us with an altogether different method for solving $D(G)$. Second, when $G$ is connected, finding $D(G,k)$ is really necessary to determine the distinguishing number of $H$ where $H = \alpha G$ (i.e., $H$ consists of $\alpha$ copies of $G$).  To distinguish $H$, each copy of $G$ must  be assigned a distinguishing labeling.  Additionally, no two copies of $G$ can be assigned equivalent distinguishing labelings.  Hence, $D(H) = \min\{k: D(G,k) \ge \alpha\}$. Finally, researchers have noted that two labels are sufficient for distinguishing many non-rigid graphs (e.g., \cite{Albertson, Tucker2, IK}).  The number of inequivalent distinguishing $k$-labelings of graphs provides one more level of granularity that enables us to differentiate between graphs with the same distinguishing numbers.    For example, consider the two graphs shown in Figure~\ref{compare-fig}.  They have the same number of vertices, their automorphism groups are isomorphic, and they can be distinguished with two labels.  Yet, $D(G_1,k) = k^4(k^4-1)/2$ but $D(G_2,k) = k^7(k-1)/2$ so with two labels at most $120$ copies of $G_1$ can be distinguished compared to $64$ copies for $G_2$. In this sense, $G_1$ is less symmetric than $G_2$ because $k$ labels can distinguish more copies of $G_1$ than $G_2$ for any $k \ge 2$. 

 To solve for $D(G,k)$, we apply two of the most common techniques for counting --  the principle of inclusion-exclusion (PIE)  and recursion.  We show that when $G$'s automorphisms are known and the size of its automorphism group, $Aut(G)$, is $O(\log n)$ where $n$ is the number of vertices in $G$ then a straightforward application of PIE can determine $D(G,k)$ efficiently.  We then modify the technique so that when $Aut(G)$ is isomorphic to $Z_t$ (the cyclic group of order $t$), $D_t$ (the dihedral group of order $2t$), $Z_t \times Z_2$, or $D_t \times Z_2$ then $D(G,k)$ can be computed in time polynomial in $n, t$ and $\log k$.  Consequently, we are able to prove that if $G$ is a triconnected planar graph then $D(G,k)$ and $D(G)$ can be determined efficiently.  
Next, by viewing $G$ via a tree decomposition $T_G$ that is made up of $G$'s cut vertices, separating pairs, and triconnected components, 
we show that $D(G,k)$ can be determined recursively. To implement this technique efficiently for a family of graphs, several ingredients are necessary including  efficient algorithms for testing graph isomorphism and finding the automorphisms of a graph's triconnected components. Since these algorithms exist for planar graphs, we arrive at the main result of the paper -- that when $G$ is a planar graph then $D(G,k)$ and $D(G)$  can be computed efficiently. 

In their introductory paper, Alberston and Collins~\cite{AC} raised the 
issue of determining the computational complexity of $DIST = 
\{(G,k)|\mbox{ $G$ has a distinguishing $k$-labeling} \}$.  
Currently, the best known result about $DIST$, which is due to  Russell and Sundaram~\cite{RS}, is that it belongs to AM, 
the set of languages for which there are Arthur and Merlin games.  
This result essentially follows from the 
fact that testing graph rigidity is in AM.  When $G$ is restricted to 
certain graph families, however, $DIST$ can belong to P.  For example, 
distinguishing numbers of cycles, hypercubes~\cite{Cowen,Chan2}, and 
acyclic graphs~\cite{AD, Cheng} can be computed efficiently.  Our main 
result extends this further -- $DIST$ belongs to P when $G$ is a planar 
graph.  Our work complements that of Fukuda, et al  
\cite{Tucker2} on triconnected planar graphs where they show that, except for  seven graphs, all graphs in this family have distinguishing number at most $2$.

 In the next section of the paper we give basic results that will be used throughout the paper.   In Section 3, we show how the principle of inclusion/exclusion  can be used to determine $D(G,k)$.  In Section 4, we develop recursive formulas for a tree decomposition of $G$ that eventually lead to the computation of $D(G,k)$. We conclude in Section 5.   We note that our algorithms  for computing $D(G,k)$ have $G$ and $k$ as input; hence, when we say that they are efficient, we mean that they run in time polynomial in the size of $G$ and $\log k$.  Additionally, these algorithms  involve addition and  multiplication.  In cases where the numbers used are functions of $k$, their values never exceed $k^n$, where $n$ is the number of nodes in graph $G$; i.e., the numbers have at most $n\log k$ bits.  Thus, in our analysis, we assume each addition takes $O(n \log k)$ time and each multiplication takes $O(n^2 \log^2 k)$ time in the worst case.

\section{Basic notions}

Suppose $\phi$ and $\phi'$ are two  distinguishing labelings of $G$.  Since (labeled) graph isomorphism is an equivalence relation, 
 we shall say that $\phi$ and $\phi'$ are {\it equivalent} if $(G,\phi) \cong (G, \phi')$; that is, there is an automorphism of $G$ that maps $(G,\phi)$ to $(G,\phi')$.  Let $\mathcal{L}(G,k)$ denote the set of all distinguishing $k$-labelings of $G$, $L(G,k)$  the size of $\mathcal{L}(G,k)$, and  $D(G,k)$ the number of equivalence classes of $\mathcal{L}(G,k)$. Below, we establish the relationships between $D(G)$, $D(G,k)$ and $L(G,k)$.

\begin{lemma}
Let $G$ be a graph and $Aut(G)$ its automorphism group. 

(i) $D(G) = \min\{k: L(G,k) > 0\} = \min\{k: D(G,k) > 0\}$.

(ii)  $D(G,k) = L(G,k)/ |Aut(G)|$.  
\label{basicproperty}
\end{lemma}

\noindent Proof:   If there is a distinguishing $k$-labeling of $G$ then the set $\mathcal{L}(G,k)$  must at least have one labeling and one equivalence class.  It follows that the smallest $k$ for which this is true must be the distinguishing number of $G$, proving the first part of the lemma. 

To prove the second part, note that $Aut(G)$ is a group that acts on $\mathcal{L}(G,k)$.  By definition, each $\phi \in \mathcal{L}(G,k)$ is preserved by only one automorphism in $Aut(G)$ -- the identity automorphism.  Hence, according to the orbit-stabilizer lemma, the size of the equivalence class of $\mathcal{L}(G,k)$ that contains $\phi$ (i.e., the orbit of $\phi$) is $|Aut(G)|$.  Consequently, the number of equivalence classes of $\mathcal{L}(G,k)$ is $L(G,k)/|Aut(G)|$. \qed

\medskip

Throughout this paper, we shall make use of Lemma~\ref{basicproperty} by viewing the problem of finding a graph's distinguishing number as a  counting problem.  While it may seem that computing $D(G,k)$ to find $D(G)$ requires more work than needed, the lemma below (first proved in \cite{Cheng}) shows that it does not if we need to distinguish multiple copies of $G$. 

\begin{lemma}
\label{countgraphs}
Let $G$ be a graph whose $g$ connected components are  $G_1, G_2, \hdots, G_g$.  Let $\phi$ be a labeling of $G$.  Then $\phi$ is distinguishing if and only if the following two conditions hold:
\begin{itemize}
\item[i.] $\phi$ when restricted to $G_i$  is distinguishing for $i=1, \hdots, g$.
\item[ii.] If $G_i \cong G_j$, $i \not = j$, then $(G_i, \phi|_{G_i}) \not \cong (G_j, \phi|_{G_j})$ for every pair of $i,j \in \{1, \hdots, g\}$.
\end{itemize}
\end{lemma}

The following is immediate.

\begin{lemma}
Let $G$ be a connected graph. If $H$ consists of  $\alpha$ copies of $G$ (i.e., $H = \alpha G$),  then $D(H) = \min \{k: D(G,k) \ge \alpha\}$.
\label{copies}
\end{lemma}

\subsection{Blocks, cut vertices, separating pairs, triconnected components}
\label{section-T_G}


Let $G=(V,E)$ be a connected graph.  Recall that $G$ is {\it $r$-connected} if  $|V| > r$ and, for any $X \subseteq V$ such that $|X| < r$, removing the vertices in $X$ from $G$ does not disconnect $G$; i.e., $G-X$ remains connected.   Suppose we are interested in determining if $G$ has some property (e.g., if it is planar).   A common technique is to first decompose $G$ into its  blocks -- which are either edges or  $2$-connected (or biconnected) subgraphs of $G$ -- and then decompose the blocks into its   ``triconnected components''\cite{HT} -- which are either parallel edges (or bonds), cycles, or $3$-connected graphs.~\footnote{Unlike blocks, however, the triconnected components of a graph need not be one of its subgraphs.}      It is then the triconnected components which are initially studied; the results are then assembled to infer the properties of the blocks, which in turn infer the property of $G$.      We shall apply this technique in  Section 4 to determine $D(G,k)$.  In particular, we shall make use of a tree, $T_G$, that captures the relationships between the cut vertices, separating pairs and triconnected components of $G$ to assemble the information for computing $D(G,k)$.


  A {\it block} of $G$ is a maximally-connected subgraph of $G$ that does not contain a cut vertex.  Thus, a block of $G$ is either an edge  or a  maximal biconnected subgraph of $G$.    Furthermore, any two blocks of $G$ have at most one vertex in common and this vertex is a cut vertex of $G$.   The {\it block-cut vertex graph} of $G$ is a bipartite graph where one partite set consists of {\it $b$-vertices} which correspond to  the blocks of $G$,  and the other partite set consists of {\it $c$-vertices} which correspond to the cut vertices of $G$.  A $b$-vertex is adjacent to a $c$-vertex if and only if  the block associated with the $b$-vertex contains  the cut vertex associated with the $c$-vertex.
  It is well known that the the block-cut vertex graph of $G$ is a tree whose leaves are all $b$-vertices and so it has  a unique center. Moreover, it can be constructed in time linear in the size of $G$~\cite{AHU}. 

Every block of $G$ that is biconnected  can similarly be represented by a tree via its triconnected components and separating pairs.  To do so,  the definition of $3$-connectedness and separating pairs have to be extended to multigraphs.  Our discussion closely follows the paper of Hopcroft and Tarjan~\cite{HT}.   Let $B$ be a biconnected multigraph, and  $\{x,y\}$ be a pair of vertices in $B$.  The set $\{x,y\}$ partitions the edge set of $B$ in the following way:  two edges belong to the same class if and only if they lie in a path that  contains neither $x$ nor $y$ except possibly as endpoints.  The classes are called the {\it separation classes of $B$ with respect to $\{x,y\}$}.  If there are at least two separation classes then the pair $\{x,y\}$ is a {\it separating pair of $B$} except when (i) there are exactly two separation classes and one class consists of a single edge, or (ii) there are exactly three classes, each consisting of a single edge.  If $B$ is a biconnected multigraph and has no separating pairs then $B$ is said to be {\it triconnected}. 

Let $\{x,y\}$ be a separating pair of $B$ and the separation classes of $B$ with respect to $\{x,y\}$ be $E_1, \hdots, E_m$.  An immediate consequence of the definition of separating pairs is that the classes can be divided into two groups  $E' = \cup_{i=1}^k E_i$ and $E''= \cup_{i=k+1}^m E_i$ so that both $E'$ and $E''$ have at least two edges.  Let $B' = ( V(E'), E' \cup \{(x,y)\})$ and $B'' = ( V(E''), E'' \cup \{(x,y)\})$. The graphs $B'$ and $B''$ are called {\it split graphs of $B$ with respect to $(x,y)$} and the edges $(x,y)$ added to both graphs are called {\it virtual edges} .  To {\it split $B$} is to replace $B$ by two of its split graphs.  Hopcroft and Tarjan suggest denoting the $i$th splitting operation via the pair $\{x,y\}$ by $s(x,y,i)$ and labeling the $(x,y)$ edges added to $B'$ and $B''$ by $i$ to differentiate this split from other splits.



 Suppose $B$ is split, its split graphs are split and so on until there are no more splits possible.  The remaining graphs are called the {\it split components of $B$}.   Clearly, they all must be triconnected; they can be grouped together as follows:  the triple bonds $\mathcal{B}_{b3}$, the (simple) triangles $\mathcal{B}_t$, and the rest of the triconnected (simple) graphs $\mathcal{B}_{tg}$.  Since there are many ways of splitting $B$, the split components of $B$ are not necessarily unique  (e.g., consider a four-cycle).   Nonetheless, this lack of uniqueness can be fixed by an operation called {\it merge} which is the reverse of split.  Let $B_1=(V_1, E_1)$ and $B_2=(V_2, E_2)$ be two split components of $B$ that contain  virtual edge $e=(x,y)$ labeled $i$. The graph $(V_1 \cup V_2, E_1 -\{e\} \cup E_2 - \{e\})$ is called the {\it merge graph of $B_1$ and $B_2$.}   To {\it merge $B_1$ and $B_2$} is to create their merge graph.  As before,  the operation is denoted by $m(x,y,i)$ to differentiate it from other merge operations.    So suppose the split components of $B$ are contained in $\mathcal{B}_{b3} \cup \mathcal{B}_t \cup \mathcal{B}_{tg}$.  Merge the triple bonds in $\mathcal{B}_{b3}$ as much as possible to  obtain a set of bonds $\mathcal{B}_b$.  Merge the triangles in $\mathcal{B}_t$ as much as possible to obtain a set of cycles $\mathcal{B}_p$.  The set of graphs in $\mathcal{B}_b \cup \mathcal{B}_p \cup \mathcal{B}_{tg}$ are called the {\it triconnected components of $B$}.  For example, a cycle has only one triconnected component -- itself -- because the triangles obtained by splitting the cycle can be merged. The following has been proven in \cite{HT}:

\begin{lemma}
\label{triconnectedcomponents}
Let $B$ be a biconnected multigraph with $m_B \ge 3$ edges. The total number of edges in the split components of $B$ is at most $3m_B - 6$.   Additionally, the  triconnected components of $B$ are unique and can be found in time linear in the size of $B$.
\end{lemma}

Lemma~\ref{triconnectedcomponents} implies that the order in which the split and merge operations are applied to decompose $B$ to its triconnected components is not important -- the same components are obtained.

 The biconnected multigraph $B$ can now be represented  by its {\it triconnected component-separating pair graph} which is a bipartite graph where one partite set consists of {\it $t$-vertices}  that correspond to the triconnected components  of $B$, 
and the other partite set consists of {\it $s$-vertices} that correspond to $B$'s separating pairs which  exist as virtual edges in $B$'s triconnected components.  A $t$-vertex is adjacent to an $s$-vertex if and only if the triconnected component associated with  the $t$-vertex contains the separating pair associated with the $s$-vertex.  It is easy to verify that this bipartite graph must again be a tree, all its leaves are $t$-vertices and consequently has a unique center.  Moreover, because the triconnected components of $B$ can be found in linear time, the tree can also be constructed in linear time. 


\begin{figure}
\label{first-example}
\centering
\epsfig{file=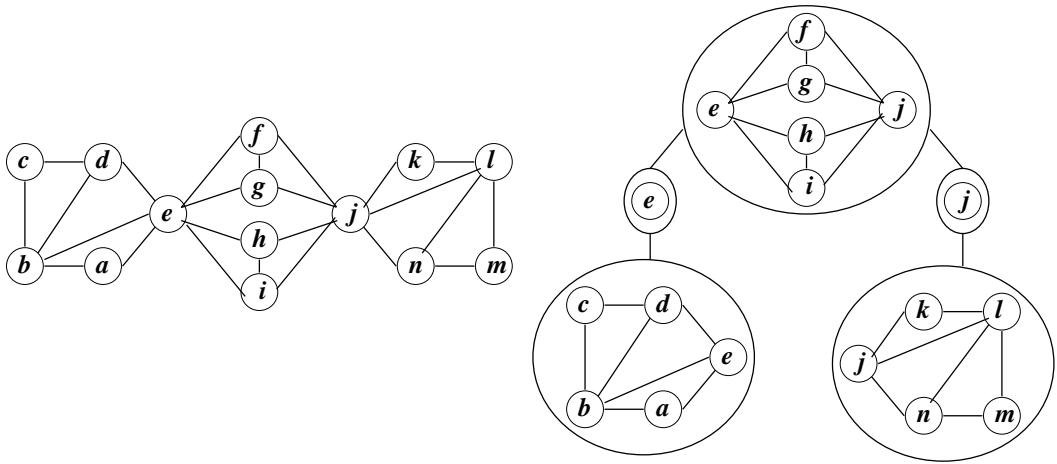, width=5.5in} 
\caption{A graph and its block-cut vertex graph.}
\end{figure}

\begin{figure}
\label{example2-fig}
\centering
\epsfig{file=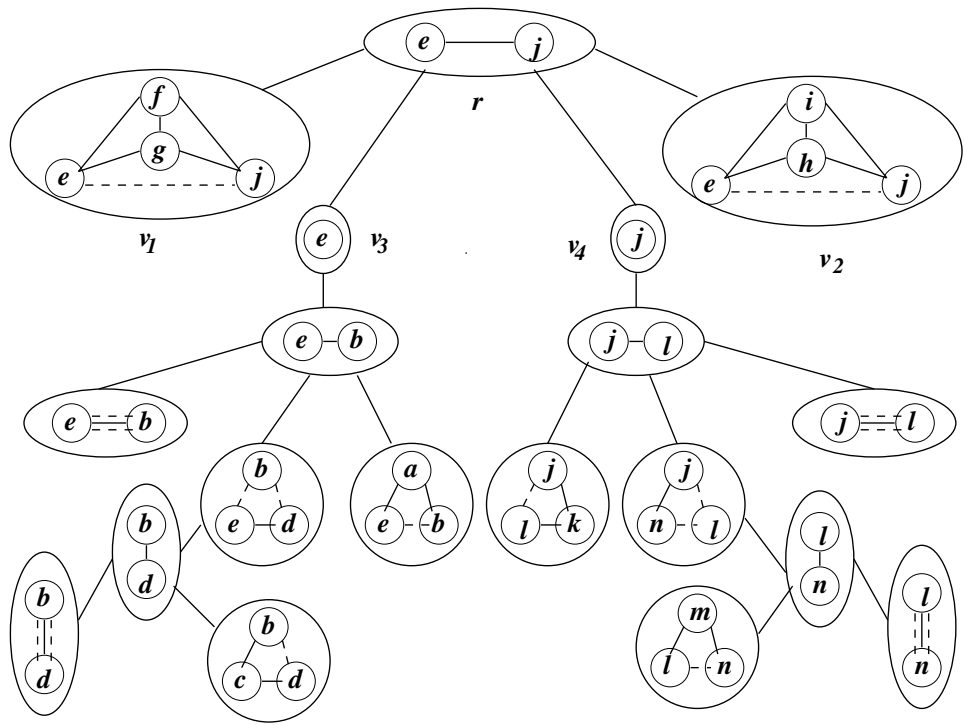, width=5in}
\caption{The tree decomposition $T_G$ of the graph in Figure 2 where $r = r(T_G)$.  The virtual edges of the triconnected components are drawn with dashed lines.}
\end{figure}

\smallskip

\noindent {\it Building a tree-decomposition of $G$.} Let $G$ be a connected graph. Let us now build a tree decomposition of $G$, $T_G$, that incorporates the triconnected component-separating pair graph of each block of $G$ into the block-cut vertex graph of $G$.  Initially set $T_G$ to be the block-cut vertex graph of $G$.  Then, for each $b$-vertex $z$ whose associated block is $B$, replace $z$ with $B$'s triconnected component-separating pair graph $T_B$.  Attach $T_B$ to each neighbor $y$ of $z$ in the following manner. Let $a$ be the cut vertex associated with $y$.   Node $a$ is part of one or more triconnected components and separating pairs of $B$. It is straightforward to check that the vertices in $T_B$ associated with these components and pairs form a subtree  which has a unique center because all the leaves of the subtree are $t$-vertices. Connect the center of this subtree in $T_B$ to $y$. 

Next, let us assign a root, $r(T_G)$, to $T_G$ as follows.  If the center of the block-cut vertex graph of $G$ is a $c$-vertex, this $c$-vertex is part of $T_G$. Set $r(T_G)$ to be this $c$-vertex.  Otherwise, the center of the block-cut vertex is a $b$-vertex associated with some block $B$.  Set $r(T_G)$ to be the center of $T_B$.  The tree decomposition of the graph in Figure 2 is shown in Figure 3.   



\begin{claim}
\label{fixroot}
Every automorphism of $G$ maps the structure associated with $r(T_G)$ -- which may be a cut vertex, a separating pair, or a triconnected component of $G$ --  to itself.  \end{claim}

\noindent Proof: Let $BC_G$ denote the block-cut vertex graph of $G$.  Recall that $BC_G$ has a unique center; denote it as $z^*$.   Every automorphism of $G$ induces an automorphism on $BC_G$.~\footnote{ That is,  if $\pi \in Aut(G)$,  define $f_{\pi}$ on the set of vertices of $BC_G$ so that $f_\pi$ mimics the actions of $\pi$ on $G$.  Thus, for each vertex $z$ in $BC_G$ whose associated structure is $A$, let  $f_{\pi}(z)$ be the vertex in $BC_G$ associated with the structure $\pi(A)$.  It is easy to verify that $f_\pi$ is an automorphism of $BC_G$.}
  But every automorphism on $BC_G$ fixes $z^*$; hence, every automorphism of $G$ fixes the structure associated with $z^*$. If $z^*$ is a $c$-vertex, $r(T_G) = z^*$ and so the claim follows.  Otherwise, $z^*$ is a $b$-vertex  that it is associated with some block $B$.  This means that the action of every automorphism of $G$ on $B$ corresponds to an automorphism of $B$.  Now, every automorphism of $B$ induces an automorphism on $T_B$. Applying the same argument above to $T_B$,  we have that every automorphism of $B$ fixes the structure associated with center of $T_B$. Since $r(T_G)$ is the center of $T_B$,  the claim follows. \qed
\smallskip

From here onwards, we shall treat $T_G$ as a rooted tree.   For each node $v$ in $T_G$,  let $T_v$  denote the subtree of $T_G$ rooted at $v$, and  $G(T_v)$  denote the graph obtained by merging (using the merge operation we defined earlier) all the triconnected components associated with the $t$-vertices in $T_v$. We make a few observations about $G(T_v)$.  When $v = r(T_G)$ then $G(T_v) = G$.  Furthermore, for general $v$, $G(T_v)$ consists of connected graphs "hanging" from the structure associated with $v$; these connected graphs can be obtained from the $G(T_w)$'s,  where $w$ is a child of $v$.  We also note  that some of the $G(T_v)$'s may not be subgraphs of $G$ -- which occurs when $v$ is an $s$-vertex and its parent is a $t$-vertex or when $v$ is a $t$-vertex and its parent is an $s$-vertex.\footnote{If $\{x,y\}$ is the separating pair associated with the $s$-vertex, then it is possible that $G(T_v)$ will contain multiple copies of the edge $(x,y)$.  We note though that we can ignore the other copies of $(x,y)$ as their multiplicity does not affect any of our computations.}





In our later discussions, we will mostly be interested in the automorphisms of $G(T_v)$ that fix a cut vertex, a separating pair, or a triconnected component, and so we use $Aut(G(T_v);*)$ to denote the set of automorphisms of $G(T_v)$ that fix the structures in $*$.  For example, let $\{x,y\}$  be a separating pair in $G(T_v)$. The automorphisms of $G(T_v)$ in $Aut(G(T_v);x,y)$ fix the vertices $x$ {\it and} $y$ while those in $Aut(G(T_v);xy)$ fix the edge $(x,y)$.  When $H$ is a triconnected component in $G(T_v)$, the automorphisms of $G(T_v)$ in $Aut(G(T_v);H)$ map $H$ to itself (i.e., the set $V(H)$ to itself) and the automorphisms in $Aut(G(T_v);H,x,y)$ map $H$ to itself and, additionally, vertices $x$ and $y$ to themselves. From Claim~\ref{fixroot}, we have the next lemma.

\begin{lemma}
\label{autroot}
Let $G$ be a connected graph and $T_G$ its tree decomposition.  Then $Aut(G) = Aut(G;A)$ where $A$ is the structure associated with $r(T_G)$.  
\end{lemma}

From the construction of $T_G$, we also have the next two lemmas.

\begin{lemma}
\label{autc-vertex}
Let $v$ be a $c$-vertex in $T_G$ and $a$ be its associated cut vertex.  Let $w$ be a child of $v$ in $T_G$.  The following are true:

\noindent (i) if $w$ is an $s$-vertex,  then it is associated with  some separating pair $\{a,b\}$  and $Aut(G(T_w); a) = Aut(G(T_w); a,b)$,

\noindent (ii) if $w$ is a $t$-vertex and its associated  triconnected component is $H$, then $H$ contains $a$ and $Aut(G(T_w);a) = Aut(G(T_w); H,a)$.  
\end{lemma}

\noindent Proof:  Since $v$ and $w$ are adjacent in $T_G$ and $v$ is a $c$-vertex while $w$ is an $s$- or a $t$-vertex, there is a block $B$ that contains cut vertex $a$ and the structure associated with $w$.  As we noted in the construction of $T_B$, $a$ must be part of one or more separating pairs and triconnected components in $B$, and  the vertices associated with these pairs and components form a subtree in $T_B$. Let us call this subtree $T_{B,a}$. Since $w$ was chosen so that it is the center of $T_{B,a}$, the structure associated with $w$ contains $a$.


By the way the block-cut vertex graph of $G$ is constructed, $B$ must be the only block in $G(T_w)$ that contains $a$.  Hence, every automorphism of $G(T_w)$ that fixes $a$ must map the separating pairs and triconnected components of $B$ that contain $a$ to similar separating pairs and triconnected components.  That is, the actions of every automorphism in $Aut(G(T_w);a)$ induces an automorphism on $T_{B,a}$.  But $T_{B,a}$ has a unique center -- $w$-- which means that every automorphism in $Aut(G(T_w);a)$ must map the structure associated with $w$ to itself.  The lemma follows.  \qed

\begin{lemma}
\label{auts-vertex}
Let $v$ be an $s$-vertex in $T_G$ and $\{x,y\}$ be its associated separating pair.  Let $w$ be a child of $v$ in $T_G$.  If $w$ is a $t$-vertex whose associated triconnected component is $H$, then $Aut(G(T_w);x,y) = Aut(G(T_w); H,x,y)$ and $Aut(G(T_w);xy) = Aut(G(T_w); H,xy)$.  
\end{lemma}

\noindent Proof:   Since $v$ and $w$ are adjacent in $T_G$ and $v$ is an $s$-vertex while $w$ is a $t$-vertex, there is again a block $B$ that contains the structures associated with both vertices.  By the way the triconnected component-separating pair graph of $B$ is constructed, it must be the case that  $H$ is the only triconnected component in $G(T_w)$ that contains $\{x,y\}$.  Hence, every automorphism in $Aut(G(T_w); x,y)$ must map $H$ to itself and so $Aut(G(T_w);x,y) = Aut(G(T_w); H,x,y)$.  By the same reasoning, $Aut(G(T_w);xy) = Aut(G(T_w); H,xy)$. \qed

\smallskip
The following lemmas will also be useful later.

\begin{lemma}
\label{T_G}
The tree $T_G$ can be constructed in $O(n^2 + nm)$ time where $n$ is the number of vertices and $m$ the number of edges in $G$.
\end{lemma}

\noindent Proof:  Constructing $G$'s block-cut vertex graph and rooting it at its center takes $O(n + m)$ time.  Creating the separating pairs-triconnected components graph $T_B$ of block $B$ takes $O(n_B + m_B)$ time where $n_B$ and $m_B$ are the number of nodes and edges in block $B$. Connecting $T_B$ to $T_G$ takes $O(c_B (n_B + m_B))$ where $c_B$ is the number of cut vertices in block $B$.  Thus, doing this for all blocks $B$ takes $O(n^2 + nm)$ time since $\sum_B n_B \leq n + m$ and $\sum_B m_B = m$. \qed
\smallskip

\begin{lemma}
\label{t-vertices-lemma}
Let $B$ be a block of $G$ with $n_B$ vertices and $m_B \ge 3$ edges.  Let $\mathcal{H}$ be the set that contains all its triconnected components of $B$.  For each $H \in \mathcal{H}$, let $\mathcal{S}_H$ denote the set containing the separating pairs of $G$ in $H$ used in the construction of $T_G$.  Then, $\sum_{H \in \mathcal{H}} |\mathcal{S}_H| = O(m_B)$ and $\sum_{H \in \mathcal{H}} |V(H)| = O(m_B)$. 
\end{lemma}

\noindent Proof:  Suppose the split operation was applied $g$ times to $B$ until no more splits are possible.  Let $\mathcal{H}'$ contain the resulting split components.   For each $H' \in \mathcal{H}'$, define $\mathcal{S}_{H'}$ as in the lemma.  We note that when $B$ is split into two components, the separating pair used to create the split becomes part of both components.  That is, each split operation contributes a value of $2$ to  $\sum_{H' \in \mathcal{H}'} |\mathcal{S}_{H'}|$.  Hence, $\sum_{H' \in \mathcal{H}'} |\mathcal{S}_{H'}| = 2g$.  Now, according to Lemma~\ref{triconnectedcomponents}, the total number of edges in the split components in $\mathcal{H}'$ is at most $3m_B -6$.  Since a split component in $\mathcal{H}'$ has at least three edges, $g \leq m_B-2$ and so $\sum_{H' \in \mathcal{H}'} |\mathcal{S}_{H'}| = O(m_B)$.  Next, notice that $|V(H')| \leq |E(H')|$ for each $H' \in \mathcal{H}'$ so $\sum_{H' \in \mathcal{H}'} |V(H')| \leq  3m_B -6$.  Finally, because $\sum_{H \in \mathcal{H}} |\mathcal{S}_H| \leq \sum_{H' \in \mathcal{H}'} |\mathcal{S}_{H'}|$ and $\sum_{H \in \mathcal{H}} |V(H)|  \leq \sum_{H' \in \mathcal{H}'} |V(H')|$, the lemma follows. \qed
\bigskip

Finally, we note that we call $T_G$ a tree decomposition of $G$ because it  really is a tree decomposition as defined by Robertson and Seymour (see Chapter 12 in \cite{Diestel} for an introduction).  That is, if $v$ is a node of $T_G$ and $V_v$ contains the vertices of the structure in $G$ associated with $v$, it should be clear from our construction that the following are true:  (i) $V_v \subseteq V(G)$ for each $v$, (ii) $\cup_v V_v = V(G)$, (iii) every edge of $G$ has two of its endpoints in some $V_v$, and (iv) whenever $y$ and $z$ are neighbors of $v$ then $V_y \cap V_z \subseteq V_v$.   In our discussion, however, it is important that we keep track of the actual structure associated with $v$ and not just the vertices in $V_v$. 

\section{Counting the distinguishing $k$-labelings of graphs via PIE}
\label{section-PIE}

Given a graph $G$ and its automorphisms, we begin by applying the principle of inclusion-exclusion (or PIE) to count  its distinguishing $k$-labelings.  Unfortunately, the technique requires the computation of $\Omega(2^{|Aut(G)|})$ terms and so becomes impractical when $G$ has many automorphisms.   We show how the method can be modified when $Aut(G)$ is isomorphic to certain groups.  In particular, we prove that when $G$ is a triconnected planar graph, $L(G,k)$, $D(G,k)$, and $D(G)$ can be computed in time polynomial in $\log k$ and the size of $G$.

Suppose $Aut(G) = \{ \pi_0, \pi_1, \hdots, \pi_{g-1} \}$ where $\pi_0$ is the identity automorphism.  Let $\phi$ be some $k$-labeling of $G$.  We say that an automorphism $\pi_i$ of $G$ {\it preserves} $\phi$ if $\phi(v) = \phi(\pi_i(v))$ for each $v$ of $G$.  Clearly, $\pi_0$ preserves $\phi$,  and if no other automorphism of $G$ preserves $\phi$ then $\phi$ is a distinguishing $k$-labeling of $G$.  Let $P \subseteq Aut(G)$ and $N_{\ge}(P)$ denote the number of $k$-labelings of $G$ that are preserved by all the automorphisms in $P$.  Let $N_=(P)$ equal the number of $k$-labelings of $G$ that are preserved by all the automorphisms in $P$ but no others. 
Thus, $L(G,k) = N_=(\{ \pi_0\})$.  According to the PIE,

\begin{eqnarray}
 N_=(\{ \pi_0 \}) & = & \sum_{ \{ \pi_0 \} \subseteq P \subseteq Aut(G)} (-1)^{|P|-1} N_\ge(P).
\label{formula1}
\end{eqnarray}

 Next,  we describe a method for computing $N_\ge(P)$, for each $P \subseteq Aut(G)$.  Suppose $\pi_i \in P$.  A $k$-labeling $\phi$ is preserved by $\pi_i$ if and only if $\phi$  assigns the same label to $v$ and to $\pi_i(v)$ for each vertex $v$ in $G$. In fact, if there is a sequence of vertices $v_1, v_2, \hdots, v_{r}$ such that $v_j = \pi_i(v_{j-1})$ for $j = 2, \hdots, r$ then $\phi$ must assign all of these $r$ vertices the same label. By extending this idea further, we arrive at the following lemma.

\begin{lemma}
Let $\pi_i \in Aut(G)$ and $\phi$ be a $k$-labeling of $G$.  Let $G_{\pi_i}$ be the graph whose node set is  $V(G)$ and whose edge set consists of the pairs $(v, \pi_i(v)), v \in V(G)$. The automorphism $\pi_i$ preserves $\phi$ if and only if, for each connected component in $G_{\pi_i}$, $\phi$ assigns the same label to all the vertices in  that component.  Consequently, let $P \subseteq Aut(G)$.  The automorphisms in  $P$ preserve $\phi$ if and only if, for each connected component in $\cup_{\pi_i \in P} G_{\pi_i}$, $\phi$ assigns the same label to all the vertices in  that component.  
\label{components}
\end{lemma}

An immediate implication of the lemma is if $\cup_{\pi_i \in P} G_{\pi_i}$ has $t$ connected components  and there are $k$ labels available then $N_\ge (P) = k^t$.  We are now ready to prove the next result.  


\begin{theorem}
\label{PIE}
Let $G$ be a graph on $n$ vertices and $k$ be a positive integer. Suppose all the automorphisms of $G$ are given.  Then  $L(G,k)$ can be computed in 
$O(n^2 \log^2k + 2^{|Aut(G)|}(n \times |Aut(G)| + n \log k))$ time.
\end{theorem}

\noindent Proof:  Begin by computing and storing the values $k, k^2, k^3, \hdots, k^n$.  Set $L(G,k)$ to $0$.  For each subset $P$ such that $\{ \pi_0 \} \subseteq P \subseteq Aut(G)$, (i) construct $\cup_{\pi_i \in P} G_{\pi_i}$ and find the number of its connected components $t$ using breadth-first-search and (ii) add $(-1)^{|P|-1} k^t$ to $L(G,k)$.  According to equation (\ref{formula1}), at the end of this algorithm the value of $L(G,k)$ is the number of distinguishing $k$-labelings of $G$.   Computing the powers of $k$ can be done in $O(n^2 \log^2 k)$ steps.  Each iteration of the for loop takes at most $O(n \times |Aut(G)| + n \log k)$ time where the first term in the sum accounts for the time it takes to construct $\cup_{\pi_i \in P} G_{\pi_i}$ and find its connected components, and the  latter term  accounts for adding $k^t$ to $L(G,k)$.  Since there are $2^{|Aut(G)|-1}$ subsets $P$ to consider, computing $L(G,k)$ takes $O(n^2 \log^2k + 2^{|Aut(G)|}(n \times |Aut(G)| + n \log k))$ time. \qed

\begin{corollary}
Let $G$ be a graph with $n$ vertices and $k$ be a positive integer. Suppose all the automorphisms of $G$ are given.   If $|Aut(G)| = O(\log n)$, then  $L(G,k)$ can be computed in time polynomial in $n$ and $\log k$. 
\end{corollary}
\medskip

The reason why implementing the PIE formula for $L(G,k)$ can take exponential time is because there are $\Omega(2^{|Aut(G)|})$ $N_\ge(P)$ terms in the formula.  Below we demonstrate that the technique can be modified when $Aut(G)$ is isomorphic to certain groups.  We consider the case when $Aut(G) \cong \Gamma$ where (i) $\Gamma = Z_t$, the cyclic group of order $t$, (ii) $\Gamma = D_t$, the dihedral group of order $2t$, and (iii) $\Gamma = Z_t \times Z_2$ or $D_t \times Z_2$.  All will be useful when we discuss triconnected planar graphs in the next subsection. A key feature of these results is that $|Aut(G)| = O(t)$ and yet the number of $N_{\ge}(P)$ terms that must be computed to derive $L(G,k)$ is polynomial in $t$, and not exponential in  $t$. Before we proceed, we first prove  the following lemma.

\begin{lemma}
Let $P \subseteq Aut(G)$ and $\langle P\rangle$ be the subgroup generated by $P$.  Every $k$-labeling of $G$ preserved by all the automorphisms in $P$ is also preserved by all the automorphisms in $\langle P \rangle$.  
\label{subgroup}
\end{lemma}

\noindent Proof: Let $\phi$ be a $k$-labeling of $G$ preserved by all the automorphisms in $P$.  Let $\pi \in \langle P \rangle$.  Since $Aut(G)$ is finite, we can write $\pi$ as $\sigma_r*\sigma_{r-1} * \hdots * \sigma_1$ where $r \in \mathbf{Z}^+$ and each $\sigma_i \in P$.  Since each $\sigma_i$ preserves $\phi$,  for each vertex $u$ of $G$, 
$$ \phi(u) = \phi(\sigma_1(u)) = \phi(\sigma_2(\sigma_1(u))) = \cdots = \phi(\sigma_r(\cdots(\sigma_2(\sigma_1(u))))).$$  That is $\pi = \sigma_r*\sigma_{r-1} * \hdots * \sigma_1$ preserves $\phi$ as well. \qed
\bigskip

In the subsequent discussion, when $Aut(G) \cong \Gamma$, we shall denote the automorphisms of $G$ as $\pi_\sigma$ where $\sigma \in \Gamma$, and let $\pi_\sigma * \pi_{\sigma'} = \pi_{\sigma*\sigma'}$. 
\medskip

\noindent {\it When $Aut(G) \cong Z_t$.}  Let $Z_t$ be the cyclic group of order $t$ and $\rho$ be one of its  generators.  Its elements are $\rho^0$ (the identity), $\rho, \rho^2, \hdots, \rho^{t-1}$ where   $\rho^i*\rho^j = \rho^{i+j \mod t}.$

\begin{theorem}
Let $Aut(G) \cong Z_t$, where the prime factorization of $t$ is $\prod_{i=1}^s p_i^{r_i}$.  Suppose a generator $\pi_{\rho}$ of $Aut(G)$ is given.   Let $P^* = \{\pi_{\rho^i} : i \in \{{t/p_1}, {t/p_2}, \hdots, {t/p_s}\} \}$.  Then 
\begin{eqnarray*}
 L(G,k) & = & \sum_{ P \subseteq   P^*} (-1)^{|P|} N_\ge(P) \label{formula2}.
\end{eqnarray*}
\label{cyclic}
\end{theorem}

\noindent Proof: To prove the theorem, we will show that a $k$-labeling $\phi$ of $G$ is distinguishing if and only if no automorphism in $P^*$ preserves $\phi$.  One direction is obvious:  if  $\phi$ is distinguishing, all non-trivial automorphisms of $G$ do not preserve $\phi$.  Since $P^*$ contains only non-trivial automorphisms of $G$, the result follows.  So suppose $\phi$ is not distinguishing.  It must be preserved by some $\pi_{\rho^j}$, $j \not = 0$.  Let $g = gcd(t,j) = \prod_{i=1}^s p_i^{t_i}$, where $0 \leq t_i \leq r_i$.  We know that  $\rho^g \in \langle \rho^j \rangle$.  Since $j < t$, we also know that $g$ must divide one of the numbers in $\{{t/p_1}, {t/p_2}, \hdots, {t/p_s}\}$, say $t/p_1$; i.e., $\rho^{t/p_1} \in \langle \rho^g \rangle$.  By Lemma~\ref{subgroup},  it follows that if  $\pi_{\rho^j}$ preserves $\phi$ then $\pi_{\rho^g}$ also preserves $\phi$, which implies that  $\pi_{\rho^{t/p_1}}$ 
does so as well.  That is, some automorphism in $P^*$ preserves $\phi$.  Applying the PIE,   $ L(G,k)  =  \sum_{P \subseteq   P^*} (-1)^{|P|} N_\ge(P).$   \qed
\smallskip


\medskip

\noindent {\it When $Aut(G) \cong D_t$.}  Let $D_t$ be the dihedral group of order $2t$.  If we let the generators of $D_t$ be the rotation $\rho$ and reflection $\tau$, then the elements of $D_t$ are $\rho^0$ (the identity), $\rho^1, \hdots, \rho^{t-1}$, $\tau \rho^0$, $\tau \rho^1, \hdots, \tau \rho^{t-1}$, where $\tau^2 = \rho^0$, $\tau \rho^i = \rho^{-i} \tau$ and $\rho^i * \rho^j = \rho^{i + j \mod t}$.

\begin{theorem}
\label{dihedral}
 Let $Aut(G) \cong D_t$, where the prime factorization of $t$ is $\prod_{i=1}^s p_i^{r_i}$.  Suppose generators $\pi_{\rho}$ and $\pi_{\tau}$ of $Aut(G)$ are given.  Let $P^* = \{\pi_{\rho^i} : i \in \{{t/p_1}, {t/p_2}, \hdots, {t/p_s}\} \}$.  Then

\begin{eqnarray*}
N_=( \{ \pi_{\rho^0},  \pi_{\tau \rho^i} \}) & = & \sum_{ \{ \pi_{\tau \rho^i} \} \subseteq P   \subseteq  \{  \pi_{\tau \rho^i} \} \cup P^*} (-1)^{|P|-1} N_\ge(P),
\label{formula3}
\end{eqnarray*}
and 
\begin{eqnarray}
 L(G,k) & = &  \sum_{P \subseteq  P^*} (-1)^{|P|} N_\ge(P) - \sum_{i=0}^{t-1} N_=( \{ \pi_{\rho^0}, \pi_{\tau \rho^i} \}).
 \label{formula4}
\end{eqnarray}

\end{theorem}

\noindent Proof: We shall first prove that a $k$-labeling $\phi$ of $G$ that is preserved by at least two non-trivial automorphisms of $G$ is also preserved by some automorphism in the set  $P^* = \{\pi_{\rho^i} : i \in \{{t/p_1}, {t/p_2}, \hdots, {t/p_s}\} \}$.   If one of the automorphisms that preserves $\phi$ is  preserved by $\pi_{\rho^j}$, $j \not = 0$, then by the proof of Theorem~\ref{cyclic} it must also be preserved by some automorphism in $P^*$.   If the two automorphisms that preserve $\phi$ are  $\pi_{\tau \rho^i}$ and $\pi_{\tau \rho^j}$, where $i < j$, then $\pi_{\tau \rho^i} * \pi_{\tau \rho^j} = \pi_{\rho^{j-i}}$ also preserves $\phi$.  Once again, some automorphism in $P^*$ must preserve $\phi$.  

To prove equation (\ref{formula4}), we now consider the set of all $k$-labelings of $G$.  Let sets  $A$, $B$, and $C$ consist of all $k$-labelings of $G$ preserved by $\pi_{\rho^0}$ only, by $\pi_{\rho^0}$ and $\pi_{\tau \rho^i}$ for some $i \in \{0, 1, \hdots, t-1\}$ only, and by some automorphism in $P^*$ respectively.  Any $k$-labeling of $G$ must belong to exactly one of the three sets because: (i)  if it is distinguishing, it belongs to set $A$ and if not to $B \cup C$; (ii)  if it is preserved by exactly one non-trivial automorphism of $G$, and it is of the form $\pi_{\tau \rho^i}$, it belongs to set $B$; otherwise, it belongs to set $C$; (iii)  finally, if it is preserved by at least two non-trivial automorphisms of $G$, then it belongs to set $C$.   That is, $A \cup B \cup C$ contains all the $k$-labelings of $G$ and no two of them have a $k$-labeling of $G$ in common.   Thus, $L(G,k) = |A| = k^n - |B| -|C|$.

By the way we defined set $B$, $|B| = \sum_{i=0}^{t-1} N_=( \{ \pi_{\rho^0}, \pi_{\tau \rho^i} \})$.  Consider a $k$-labeling of $G$ that is preserved by $\pi_{\tau \rho^i}$.  From our earlier argument, we can assume that such a $k$-labeling is preserved by $\pi_{\tau \rho^i}$ only or by $\pi_{\tau \rho^i}$ and some other automorphism in $P^*$, in addition to being preserved by $\pi_{\rho^0}$.   According to the PIE, this means that  $N_=(\{  \pi_{\rho^0},\pi_{\tau \rho^i} \})  =  \sum_{ \{  \pi_{\tau \rho^i} \} \subseteq P   \subseteq  \{  \pi_{\tau \rho^i} \} \cup P^*} (-1)^{|P|-1} N_\ge(P)$.  Finally, $C$ consists of all the $k$-labelings of $G$ preserved by at least one of the automorphisms in $P^*$.  So, according to the PIE,  $ k^n - |C|   = \sum_{P \subseteq P^*} (-1)^{|P|}N_\ge(P)$.   Hence, $L(G,k) = \sum_{P \subseteq P^*} (-1)^{|P|}N_\ge(P) - \sum_{i=0}^{t-1} N_=( \{  \pi_{\rho^0}, \pi_{\tau \rho^i} \})$, which proves equation (\ref{formula4}).  \qed

\medskip

\noindent {\it Example.} Consider the cycle on $n$ vertices $C_n$ where $n$ is a prime number.  Then $Aut(C_n) = D_{n}$ and $P^* = \{ \pi_{\rho} \}$.  To solve for $L(C_n,k)$, we need the following values:  $N_\ge(\emptyset)$, $N_\ge(\pi_{\rho})$, $N_\ge(\{ \pi_{\tau \rho^i}, \pi_{\rho} \})$ and $N_\ge(\{ \pi_{\tau \rho^i} \})$.  Every $k$-labeling of $C_n$ should be counted in $N_\ge(\emptyset)$ so $N_\ge(\emptyset) = k^n$. To solve for $N_\ge(\pi_{\rho})$, recall that we considered $G_{\pi_{\rho}}$ which is a graph that has only one component.  Hence, $N_\ge(\pi_{\rho}) = k$.  Similarly, $N_\ge(\{ \pi_{\tau \rho^i}, \pi_{\rho} \}) = k$.  Finally, $G_{ \pi_{\tau \rho^i} }$ consists of $(n+1)/2$ components  since any reflection of $C_n$ fixes one vertex $v$ and maps the equidistant vertices from $v$ to each other.  Thus, $N_\ge(\{ \pi_{\tau \rho^i} \}) = k^{(n+1)/2}$.  From  equation (\ref{formula4}), 
\begin{eqnarray*}
L(C_n,k) & = &  N_\ge(\emptyset) - N_\ge(\{\pi_{\rho}\}) - \sum_{i=0}^{n-1}\left( N_\ge(  \{ \pi_{\tau \rho^i} \}) - N_\ge(\{ \pi_{\tau \rho^i}, \pi_{\rho} \}) \right) \\
      & = &  k^n - k - nk^{(n+1)/2} + nk \\
      & = &  k^n - nk^{(n+1)/2} + (n-1) k \\
      & = &  k (k^{(n-1)/2} - 1) (k^{(n-1)/2} - (n-1)).
\end{eqnarray*}

Consequently,  $D(C_n,k) = k (k^{(n-1)/2} - 1) (k^{(n-1)/2} - (n-1)) / 2n$. When $n = 5$, for example, $D(C_5,1) = D(C_5, 2) = 0$ but $D(C_5,3) = 12$ so $D(C_5) = 3$.

\medskip

\noindent {\it When $Aut(G) \cong Z_t \times Z_2$ or $D_t \times Z_2$.}  We state the following theorem without proof because the arguments are just extensions of those in Theorems~\ref{cyclic} and \ref{dihedral}.

\begin{theorem}
\label{crossproduct}
 Suppose the prime factorization of $t$ is $\prod_{i=1}^s p_i^{r_i}$, the group $Z_t \times Z_2 = \{(\rho^i, \sigma^j), i \in \{0,1, \hdots t-1\}, j \in \{0,1\} \}$ and the group $D_t \times Z_2 = \{(\rho^i, \sigma^j), (\tau\rho^i, \sigma^j), i \in \{0,1, \hdots m-1\}, j \in \{0,1\} \}$.  When $t$ is odd, set $P^*_0 = \{ \pi_{(\rho^0, \sigma)} \}$; otherwise, set $P^*_0 = \{ \pi_{(\rho^0, \sigma)}, \pi_{(\rho^{t/2}, \sigma)} \}$.  Let $P^* = P^*_0 \cup \{ \pi_{(\rho^i, \sigma^0)}: i \in \{ t/p_1, t/p_2, \hdots, t/p_s  \} \}$. \\

\noindent (i) When $Aut(G) \cong Z_t \times Z_2$, 
\begin{eqnarray*}
L(G,k) & = & \sum_{P \subseteq P^*} (-1)^{|P|} N_\ge (P).
\end{eqnarray*}

\noindent (ii)  When  $Aut(G) \cong D_t \times Z_2$, and for $b=0$ or $1$,

\begin{eqnarray*}
N_=(\{  \pi_{(\rho^0, \sigma^0)}, \pi_{(\tau \rho^i, \sigma^b)} \}) & = & \sum_{ \{ \pi_{(\tau \rho^i, \sigma^b)} \} \subseteq P \subseteq  \{ \pi_{(\tau \rho^i, \sigma^b)} \}    \cup P^* } (-1)^{|P|-1} N_\ge(P)
\end{eqnarray*}
and 

\begin{eqnarray*}
L(G,k) & = & \sum_{P \subseteq P^*} (-1)^{|P|} N_\ge(P) - \sum_{b=0}^1 \sum_{i=0}^{t-1} N_=(\{ \pi_{(\rho^0, \sigma^0)},\pi_{(\tau \rho^i, \sigma^b)} \}).
\end{eqnarray*}

\end{theorem}
\medskip

\noindent {\it Remark:} Since the number of prime factors of $t$ is 
$O(\log t)$, the number of $N_{\ge}(P)$ terms in the formula for computing 
$L(G,k)$ is $O(t)$ when $Aut(G) \cong Z_t$ or $Z_t \times Z_2$, and 
$O(t^2)$ when $Aut(G) \cong D_t$ or $D_t \times Z_2$.

\subsection{When $G$ is a triconnected planar graph}

 What is interesting about the family of triconnected planar graphs is that the automorphism groups of the graphs are only of limited kinds.  


\begin{fact}
\label{fact-lemma} \cite{Mani}
Let $G$ be a triconnected planar graph. The automorphism group of $G$ is isomorphic to a subgroup  of one of the following groups: $A_4$, $A_5$, $S_4$, $A_4 \times Z_2$, $A_5 \times Z_2$, $S_4 \times Z_2$,  $Z_t$, $D_{t}$, $Z_t \times Z_2$, $D_{t} \times Z_2$, for some integer $t$. 
\end{fact}

Since a subgroup of a dihedral group is a cyclic group or a dihedral group, clearly the subgroups of $D_t \times Z_2$ are cylic, dihedral or isomorphic to $Z_{t'} \times Z_2$ or $D_{t'} \times Z_2$ where $t' \leq t$.  In other words, the automorphism group of a triconnected planar graph is either bounded by a constant or it is isomorphic to one of four groups only.

Additionally, because triconnected planar graphs have only unique embeddings on the plane up to equivalence~\footnote{A triconnected planar graph can have two planar embeddings one of which is a mirror image of the other.}, finding all their automorphisms can also be done efficiently.  We sketch one such  method next.   Let $G$ be a triconnected planar graph with $n$ vertices and $m$ edges.  Let $e=(u,v)$ be an edge of $G$.  Let us designate its direction as being from $u$ to $v$ and one of the faces $F$ that it borders as its right face.  Create a copy of $G$, $G_{e,F}$, which specially marks $e$ and its direction, and face $F$.  For any edge $e'=(u',v')$ whose direction and right face $F'$ is fixed, create an analogous graph $G_{e',F'}$, and  using a planar graph  isomorphism testing algorithm determine if $G_{e,F}$ and $G_{e',F'}$ are isomorphic (where the marked edge and face of $G_{e,F}$ are mapped to the marked edge and face of $G_{e',F'}$).  If so, then there is an automorphism of $G$ that maps $e$ to $e'$ and $F$ to $F'$; moreover, by visiting the faces of $G_{e,F}$ and $G_{e',F'}$ in the same order, the rest of $\pi$ can be determined in  time linear in the size of $G$.   Since there is a linear time isomorphism testing algorithm for planar graphs~\cite{HW}, each iteration of the for loop takes $O(n)$ time. And since there are $O(m)$ iterations then in $O(nm) = O(n^2)$ time all the automorphisms of $G$ can be determined.  Furthermore, because each edge has two directions and two faces bordering it, the algorithm above also shows that  $|Aut(G)| \leq 4m = O(n)$ when $G$ is a triconnected planar graph.   
\medskip

To solve for $L(G,k)$ for triconnected planar graphs, we do the following: if $|Aut(G)| \leq 5! $,  use Theorem~\ref{PIE}.   Otherwise, determine if $Aut(G)$ is cyclic, dihedral, isomorphic to a direct product of a cyclic group and $Z_2$, or to a direct product of a dihedral group and $Z_2$.   If $Aut(G)$ is cyclic or dihedral,  apply Theorems~\ref{cyclic} or \ref{dihedral} respectively;  otherwise, apply Theorem~\ref{crossproduct}.

\begin{theorem}
\label{triconnected-theorem}
Let $G$ be an $n$-vertex triconnected planar graph.  Computing $L(G,k)$ and  $D(G,k)$ can be done in $O(n^2 \log^2 k + n^3 \log n + n^3 \log k)$ time.   Consequently, computing $D(G)$ takes $O(n^3 \log^2 n)$ time.  
\end{theorem}

\noindent Proof:  As we stated earlier, if $G$ has at most $5!$ automorphisms, we  use Theorem~\ref{PIE} to solve for $L(G,k)$ and $D(G,k)$.  Otherwise, we need to determine which of the four groups $Aut(G)$ is isomorphic to. In particular, $Aut(G)$ falls into CASE $i$  where $i = 1$ if the group is cyclic, $i=2$ if the group is isomorphic to $Z_t \times Z_2$ for some $t$,  $i=3$ if the group is dihedral, and $i=4$ if the group is isomorphic to $D_t \times Z_2$ for some $t$. We note that there is some  overlap in the four cases because if $t$ is odd, $Z_t \times Z_2 \cong Z_{2t}$ and $ D_{2t} \times Z_2  \cong D_{4t}$.  Thus, when we say that $Aut(G)$ belongs to CASE $2$ or $4$, we shall assume that $t$ is even.  We describe our algorithm $TriconnectCount(G,k)$ in Figure~\ref{triconnectalg}. 

 In the first part of our algorithm, we determine the case which $Aut(G)$ belongs to by considering the order of each element in $Aut(G)$.  It is easy to verify the following facts:  (i) if $Aut(G)$ has an element with order $|Aut(G)|$ it must be cyclic, (ii) if $Aut(G)$ has only three elements with order $2$ (and $3 < |Aut(G)|/2$) then it belongs to case $2$, (iii) if $Aut(G)$ has between $|Aut(G)|/2$ and $|Aut(G)|/2 + 1$ of its elements with order $2$, it  belongs to case $3$.  Once the appropriate case for $Aut(G)$ is determined, we set the value of $t$.  

The second part of the algorithm begins by computing the prime factors of $t$, finding an element $\pi \in Aut(G)$ such that  the order of $\pi$ is $t$, and then computing $P^* =  \{ \pi^i: i \in \{t/p_1, t/p_2, \hdots, t/ p_s \} \}$.  If $Aut(G)$ is cyclic  or dihedral, $P^*$ is indeed the  one needed in Theorems~\ref{cyclic} and \ref{dihedral} respectively to compute $L(G,k)$.  In cases $2$ and $4$, two more elements are missing in $P^*$.  To understand what they are, we note that  since $t$ is even $\pi$ would be of the form $\pi_{(\rho, \sigma^b)}$ where $b=0$ or $1$, and $\rho$ and $\sigma$ are generators of $Z_t$ and $Z_2$ respectively.  If we set $p_1 = 2$, then  $(\pi_{(\rho, \sigma^b)})^{t/2} = \pi_{(\rho^{t/2}, \sigma^0)}$  or $\pi_{(\rho^{t/2}, \sigma^1)}$, and $(\pi_{(\rho, \sigma^b)})^{t/p_i} = \pi_{(\rho^{t/p_i}, \sigma^0)}$   for $i = 2, \hdots s$.  At this point, the two missing elements in $P^*$ have order $2$; they can be distinguished from the other elements of $Aut(G)$ with order $2$ because they commute with every other element of $Aut(G)$ (i.e., they belong to the center of $Aut(G)$), whereas the others do not.  By updating $P^*$, we now obtain the appropriate $P^*$ in Theorem~\ref{crossproduct}.   Finally, for cases $3$ and $4$, we place all elements of $Aut(G)$ with order $2$ not in $P^*$ into set $T$.  It is easy to check that the rest of the algorithm computes $L(G,k)$ correctly since they follow directly from the theorems we have established. 

Computing and storing the powers of $k$ takes $O(n^2 \log^2k)$ time.  Finding all the automorphisms of $G$ take $O(n^2)$ time.  It is easy to verify that in the rest of the algorithm, the bottleneck is in computing the value of $L(G,k)$ when $|Aut(G)| > 5!$.  Applying the same analysis we used in Theorem~\ref{PIE}, and noting that $|P^*| = O(\log t)$ and $|T| = O(t)$, computing $L(G,k)$ takes $O(t^2(n \log t + n \log k))$ time.  Finally, because $G$ is a triconnected graph $|Aut(G)| = O(n)$ so $t = O(n)$.  Hence, the total runtime of  $TriconnectCount(G,k)$ is $O(n^2 \log^2 k + n^3 \log n + n^3 \log k)$.  Once  we have the value for $L(G,k)$, we also know $D(G,k)$.  To find $D(G)$, do a binary search over the range $[1,n]$ to determine the smallest $k$ for which $D(G,k) > 0$ to find $D(G)$.  The runtime in the theorem follows.  \qed

\begin{figure}

\noindent $TriconnectCount(G, k)$ \\
\noindent Input:  A triconnected planar graph $G$ with $n$ vertices, a positive integer $k$. \\
\noindent Output:  The value of $L(G,k)$.
\medskip

\noindent Compute and store the values $k, k^2, k^3, \hdots, k^n$.

\noindent Find all the automorphisms of $G$.

\noindent If $|Aut(G)| \leq 5!$ \\
\hspace*{.2in} $L(G,k) \leftarrow  \sum_{\{ \pi_0 \} \subseteq P \subseteq Aut(G)} (-1)^{|P|-1} N_\ge(P)$\\
\hspace*{.2in} return($L(G,k)$)\\
\noindent else  \\
\hspace*{.2in} compute the order of each automorphism $\pi \in Aut(G)$  \\
\hspace*{.2in} if there is an automorphism whose order is $|Aut(G)|$ \\
 \hspace*{.4in} $\mbox{CASE} \leftarrow 1$, $t \leftarrow |Aut(G)|$,  \\
 \hspace*{.2in} else \\
 \hspace*{.4in}if there are only $3$ automorphisms with order $2$ \\
 \hspace*{.6in} $\mbox{CASE} \leftarrow 2$, $t \leftarrow |Aut(G)|/2$,  \\
 \hspace*{.4in} else \\
 \hspace*{.6in} if there are between $|Aut(G)|/2$ and $|Aut(G)|/2 + 1$ elements with order $2$\\  \hspace*{.8in}  $\mbox{CASE} \leftarrow 3$, $t \leftarrow |Aut(G)|/2$,  \\
 \hspace*{.6 in} else \\
 \hspace*{.8in} $\mbox{CASE} \leftarrow 4$, $t \leftarrow |Aut(G)|/4$.  \\
 \hspace*{.2in} Compute the prime factors of $t$: $p_1, p_2, \hdots, p_s$. \\  
 \hspace*{.2in} Find an automorphism $\pi \in Aut(G)$ whose order is $t$.  \\
 \hspace*{.2in} Compute $P^* = \{ \pi^i: i \in \{t/p_1, t/p_2, \hdots, t/ p_s \} \}$.   \\
 \hspace*{.2in}  If $ \mbox{CASE} = 2$ or $4$ \\
 \hspace*{.4in}  add to $P^*$ the two automorphisms of $G$ which belong to the center of $Aut(G)$ not yet in $P^*$. \\
  \hspace*{.2in}  If $ \mbox{CASE} = 3$ or $4$ \\
 \hspace*{.4in} let $T$ consist of all automorphisms in $Aut(G)$ that is not in $P^*$ whose order is $2$. \\
 \hspace*{.2in}  $L(G,k) \leftarrow \sum_{P \subseteq P^*} (-1)^{|P|} N_\ge(P)$.   \\
 \hspace*{.2in} If $ \mbox{CASE} = 1$ or $2$ \\
 \hspace*{.4in} return($L(G,k)$)\\
 \hspace*{.2in} else \\
 \hspace*{.4in} while $T \not = \emptyset$ \\
 \hspace*{.6in}pick $\pi' \in T$ and delete $\pi'$ from $T$  \\
 \hspace*{.6in} $L(G,k) \leftarrow L(G,k) - \sum_{\{\pi'\} \subseteq P \subseteq \{\pi'\} \cup P^*} (-1)^{|P|-1} N_\ge(P)$ \\
 \hspace*{.4in} return ($L(G,k)$).   

 \label{triconnectalg}

 \caption{The algorithm for computing the number of distinguishing $k$-labelings of a triconnected planar graph.}
 \end{figure}

\section{Computing $D(G,k)$ via recursion}


In this section, we shall generalize the recursive technique (discovered independently by Arvind and Devanur\cite{AD} and by Cheng\cite{Cheng}) that was used to compute the distinguishing numbers of trees.   The main idea behind the technique is quite simple.  Let $T$ be a tree rooted at $r$.  Let $T_v$ denote the subtree of $T$ rooted at vertex $v$.  Start by setting $D(T_v,k) = k$ for each leaf $v$ since a single node has $k$ distinguishing $k$-labelings. Then, for $i = \mbox{depth$(T) - 1$ to $0$}$,  do the following:  for all nodes $v$ at depth $i$, compute $D(T_v,k)$ based on the values computed for $D(T_w,k)$, $w$ a child of $v$ in $T$.  Thus, at the end of the algorithm $D(T_r,k)$,  which equals $D(T,k)$, is determined.   To apply the above technique to a connected graph $G$, we will view $G$ as rooted tree using the tree decomposition $T_G$ described in Section~\ref{section-T_G}.  Additionally, we will also consider a generalized version of the distinguishing $k$-labelings of a graph which we shall define shortly.  Finally, we need to develop recursive formulas that relate the number of (generalized)  inequivalent distinguishing $k$-labelings of $G(T_v)$ with those of $G(T_w)$, $w$ a child of $v$ in $T_G$.

 Let $\Gamma$ be a subgroup of $Aut(G)$.   We say that a labeling $\phi$ of $G$ is {\it $\Gamma$-distinguishing} if no non-trivial automorphism in $\Gamma$ preserves $\phi$, and that two labelings $\phi$ and $\phi'$ of $G$ are {\it equivalent with respect to $\Gamma$} if some automorphism in $\Gamma$ maps $(G,\phi)$ to $(G, \phi')$.    Let $\mathcal{L}(G,k;\Gamma)$ be the set consisting of the $\Gamma$-distinguishing $k$-labelings of $G$, $L(G,k;\Gamma)$ be the size of $\mathcal{L}(G,k;\Gamma)$, and $D(G,k;\Gamma)$ be the number of equivalence classes of $\mathcal{L}(G,k;\Gamma)$ with respect to $\Gamma$.  
 When $\Gamma = Aut(G;*)$ as defined in Section~\ref{section-T_G}, we shall refer to $\mathcal{L}(G,k;\Gamma)$, $L(G,k;\Gamma)$ and $D(G,k;\Gamma)$ as $\mathcal{L}(G,k;*)$, $L(G,k;*)$, and $D(G,k;*)$ respectively.   Finally, when $(x,y)$ is an edge of $G$, we will at times differentiate between the case when a $k$-labeling of $G$ assigns $x$ and $y$ the same or different colors.  When we do so, we will place a subscript next to $\mathcal{L}$, $L$, and $D$; the subscript is $1$ if $x$ and $y$ are assigned the same color and is $2$ otherwise.  Thus, $\mathcal{L}_1(G,k;xy)$ consists of all $k$-labelings of $G$ in $\mathcal{L}(G,k;xy)$ that assigned $x$ and $y$ the same color, etc.   It is easy to verify that the following version of Lemma~\ref{basicproperty} remains true:

\begin{lemma}
\label{basicproperty2}
Let $G$ be a graph and $\Gamma$ be a subgroup of $Aut(G)$.  Then $D(G,k;\Gamma) = L(G,k;\Gamma)/|\Gamma|$.
\end{lemma}

Given a connected graph $G$,  we showed in Section 2.1 how to construct a tree decomposition of $G$, $T_G$. The construction started with $G$'s block-cut vertex graph.  Each $b$-vertex whose associated block is $B$ is then replaced with $B$'s triconnected component-separating pair graph $T_B$ and then connected to the rest of block-cut vertex graph.  Thus, $T_G$ is made up of $c$-, $s$-, and $t$-vertices which represent the cut vertices, separating pairs and triconnected components of $G$.  We shall now describe recursive formulas for $D(G(T_v),k; *)$ based on the type of vertex $v$ is in $T_G$.


\begin{theorem}
\label{c-vertex}
Let $v$ be a $c$-vertex in $T_G$ and $a$ be the cut vertex in $G$ associated with $v$.  Suppose when all the graphs in $\mathcal{G} = \{G(T_w): w \mbox{ is a child of $v$ in $T_G$} \}$ are fixed at $a$, there are  $g$ isomorphic classes and the $i$th isomorphic class contains $m_i$ copies of the connected graph $G_i$; i.e., $\mathcal{G} = m_1 G_1 \cup m_2 G_2 \cup \hdots \cup m_g G_g$. Then 
$$  D(G(T_v),k; a) = k \prod_{i=1}^g \binom{D(G_i,k; a)/k}{m_i}. $$
\end{theorem}

\noindent Proof:  By the way  $T_G$ was constructed, if $\mathcal{G} = m_1 G_1 \cup m_2 G_2 \cup \hdots \cup m_g G_g$ then $G(T_v)$ is made up of $\sum_{i=1}^g m_i$ connected components all hanging from vertex $a$.  It is easy to verify that $\phi$ is a labeling in $\mathcal{L}(G(T_v),k;a)$ if and only if $\phi$ assigns inequivalent labelings from $\mathcal{L}(G_i,k;a)$ to the $m_i$ copies of $G_i$, $i = 1, \hdots, g$ and the labels assigned to vertex $a$ by all  the labelings are the same.  This means that an equivalence class of $\mathcal{L}(G(T_v),k;a)$ is defined by (i) the label assigned to $a$ and (ii) the set of $m_i$ equivalence classes from $\mathcal{L}(G_i,k;a)$ that contain the labelings of the $m_i$ copies of $G_i$, $i = 1, \hdots, g$.  There are $k$ possible labels for $a$. Once the label for $a$ is chosen say $l$, there are $D(G_i,k; a)/k$ different equivalence classes of $\mathcal{L}(G_i,k;a)$ which assign vertex $a$ the same label. This is so because the number of equivalence classes of $\mathcal{L}(G_i,k;a)$ where $a$ is assigned the label $l$  must be the same for every possible value of $l$. It follows that there are $\binom{D(G_i,k;a)/k}{m_i}$ different sets of $m_i$ equivalence classes  that can contain the labelings assigned to the $m_i$ copies of $G_i$ for $i = 1, \hdots, g$. By the product rule of counting, the theorem is established.   \qed
\medskip

The next two theorems deal with the case when $v$ is an $s$-vertex.

\begin{theorem}
\label{s-vertex-1}  
Let $v$ be an $s$-vertex in $T_G$ and $\{x,y\}$ be the separating pair associated with $v$.  If it exists, let $w_x$ denote the child of $v$ that is a $c$-vertex associated with $x$.  Similarly, if it exists, let $w_y$ denote the child of $v$ that is a $c$-vertex associated with $y$.  Suppose when all the graphs in $\mathcal{G} = \{ G(T_w):w \mbox{ is a $t$-vertex and a child of $v$ in $T_G$} \}$ are fixed at $x$ and $y$,  there are $g$ isomorphic classes and the $i$th isomorphic class has $m_i$ copies of the connected graph $G_i$; i.e. $\mathcal{G} = m_1 G_1 \cup m_2 G_2 \cup \hdots \cup m_g G_g$. Then $D(G(T_v),k; x,y)$ equals
\begin{eqnarray*}
 &  & k^2 \max\{D(G(T_{w_x}),k;x)/k,1 \}  \max\{D(G(T_{w_y}),k;y)/k,1 \}\prod_{i=1}^g \binom{D(G_i,k; x,y)/k^2}{m_i} .
\label{fixxfixy}
\end{eqnarray*}
\end{theorem}

\noindent Proof:  Once again, it is easy to verify that $\phi \in \mathcal{L}(G(T_v),k;x,y)$ if and only if $\phi$ assigns inequivalent labelings from $\mathcal{L}(G_i,k;x,y)$ to the $m_i$ copies of $G_i$, $i=1, \hdots, g$ and the labels assigned to $x$ and to $y$ by all the labelings are the same.  Thus,
 an equivalence class of $\mathcal{L}(G(T_v),k;x,y)$ is defined by (i) the labels assigned to $x$ and $y$,  (ii) the equivalence classes of $\mathcal{L}(G(T_{w_x}),k;x)$ and $\mathcal{L}(G(T_{w_y}),k;y)$ that contain the labelings of $G(T_{w_x})$ and $G(T_{w_y})$ respectively, and (iii) the set of $m_i$ equivalence classes of $\mathcal{L}(G_i,k;x,y)$ that contain the labelings of the $m_i$ copies of $G_i$, $i = 1 \hdots g$.     There are $k^2$ labels available for $x$ and $y$.  Once the labels are chosen say $l_x$ and $l_y$, the number of equivalence classes of $\mathcal{L}(G_i,k;x,y)$ where $x$ and $y$ are assigned the said labels is $D(G_i,k; x,y)/k^2$. This is so because  in any labeling in $\mathcal{L}(G_i,k;x,y)$,  the labels of the vertices in $G_i$ other than $x$ and $y$ are the ones that actually destroy the non-trivial automorphisms of $Aut(G_i; x,y)$. Consequently, the number of labelings and the number of equivalence classes of $\mathcal{L}(G_i,k;x,y)$  with $x$ and $y$ assigned $l_x$ and $l_y$ must be  the same regardless of the pair $(l_x, l_y)$.   Of the $D(G_i,k; x,y)/k^2$ equivalence classes of $\mathcal{L}(G_i,k; x, y)$ that are being considered, $m_i$ must be chosen to contain the labelings of the $m_i$ copies of $G_i$. 
 Similarly, the number of equivalence classes of $\mathcal{L}(G(T_{w_x}),k;x)$ where $x$ is assigned the label $l_x$ and the number of equivalence classes of $\mathcal{L}(G(T_{w_y}),k;y)$ where $y$ is assigned the label $l_y$ are  $D(G(T_{w_x}),k;x)/k$ and $D(G(T_{w_y}),k;y)/k$ respectively.  By the product rule of counting, the theorem is established. \qed
\smallskip

\begin{figure}
\label{counterexample}
\centering
\epsfig{file=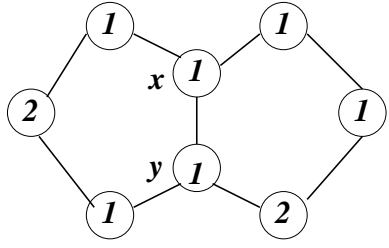, width=2.0in}
\caption{In this graph, $(x,y)$ is a separating pair.  Notice that the labeling of the $5$-cycle on the left does not destroy the automorphism of the $5$-cycle that flips the graph along the edge $(x,y)$ but the labeling for the entire graph is distinguishing.}
\end{figure}

 We will also need to compute $D(G(T_v),k;xy)$. Unlike our previous characterizations, however, it is not necessarily the case that when $\phi \in \mathcal{L}(G(T_v),k;xy)$ then $\phi$ assigned inequivalent labelings from $\mathcal{L}(G_i,k;xy)$ to each copy of $G_i$.  Figure~\ref{counterexample} shows one such exception.  Our approach this time is to consider the equivalence classes of $\mathcal{L}(G(T_v),k;x,y)$ and count those that {\it do not} belong to $\mathcal{L}(G(T_v),k;xy)$.  

Consider an arbitrary graph $H$ with an edge $(x,y)$ and suppose 
$Aut(H;xy) \not = Aut(H;x,y)$. Let $Aut(H; x\rightarrow y, y\rightarrow 
x)$ denote the set of automorphisms of $H$ that map $x$ to $y$ and $y$ to 
$x$.  Notice that $Aut(H;xy)$ is the disjoint union of $Aut(H;x,y)$ and 
$Aut(H;x\rightarrow y, y \rightarrow x)$.  Moreover, $Aut(H;x,y)$ is a 
subgroup of $Aut(H;xy)$ and $Aut(H; x\rightarrow y, y \rightarrow x)$ is a 
coset of $Aut(H;x,y)$.


Next, consider the equivalence classes of $\mathcal{L}(H,k;x,y)$.  For 
each class, either all the labelings belong to $\mathcal{L}(H,k;xy)$ 
(i.e., they destroy all the automorphisms in $Aut(H;xy)$)  or all do not.  Let 
$\mathcal{B}(H,k;x,y)$ be the set that contains all equivalence classes of
$\mathcal{L}(H,k;x,y)$ whose labelings do not belong to 
$\mathcal{L}(H,k;xy)$, and denote its size as $B(H,k;x,y)$. Our discussion will focus on computing $B(H,k;x,y)$ because this is the number of "bad" equivalence classes of  $\mathcal{L}(H,k;x,y)$ in that they do not carry over as equivalence classes of $\mathcal{L}(H,k;xy)$. Suppose $\pi 
\in Aut(H; x\rightarrow y, y \rightarrow x)$.  When all the labelings in 
an equivalence class of $\mathcal{L}(H,k;x,y)$ destroy all the 
automorphisms in $Aut(H;xy)$ then $\pi$ maps these labelings to the 
labelings of another equivalence class of $\mathcal{L}(H,k;x,y)$.  On the 
other hand, when all the labelings in an equivalence class of 
$\mathcal{L}(H,k;x,y)$ do not destroy the automorphisms in $Aut(H;xy)$ 
then $\pi$ maps these labelings to themselves.  In other words, under the 
action of $Aut(H;xy)$ the equivalence classes of $\mathcal{L}(H,k;x,y)$ 
either get paired up or stay singleton.  The ones that get paired up are 
precisely the equivalence classes of $\mathcal{L}(H,k;xy)$; i.e., each 
equivalence class of $\mathcal{L}(H,k;xy)$ is made up of two equivalence 
classes of $\mathcal{L}(H,k;x,y)$.  We shall say that such a pair of 
equivalence classes are {\it partners} in $\mathcal{L}(H,k;x,y)$.  The 
ones the stay single are the equivalence classes in $\mathcal{B}(H,k;x,y)$. 
We have proved the following lemma.

\begin{lemma}
\label{s-vertex-lemma-1}
Let $H$ be a graph with edge $(x,y)$ and $Aut(H;xy) \not = Aut(H;x,y)$.  Then $D(H,k;x,y) = 2D(H,k;xy) + B(H,k;x,y)$.
\end{lemma}

We also need the following lemma.

\begin{lemma}
\label{s-vertex-lemma-2}
Let $H$ be a graph with edge $(x,y)$  and $Aut(H;xy) \not = Aut(H;x,y)$.  Then $D_2(H,k;xy) = \binom{k}{2} D(H,k;x,y)/k^2 $ so   $D_1(H,k;xy) = D(H,k;xy) - (k-1)D(H,k;x,y)/2k$.
\end{lemma}

\noindent Proof: Consider the labelings in $\mathcal{L}_2(H,k;xy)$.  Since $x$ and $y$ are assigned different labels, they immediately destroy all the automorphisms in $Aut(H; x\rightarrow y, y\rightarrow x)$.  Thus, we can construct a labeling in $\mathcal{L}_2(H,k;xy)$ by first choosing distinct labels for $x$ and $y$, then choosing the equivalence classes of $\mathcal{L}(H,k;x,y)$ that will contain the labeling of $H$, and finally picking the labeling of $H$ from the equivalence class.  Hence, $L_2(H,k;xy) = k(k-1)\times D(H,k;x,y)/k^2 \times |Aut(H,k;x,y)|$ so 
\begin{eqnarray*}
D_2(H,k;xy) & = & L_2(H,k;xy)/ |Aut(H,k;xy)| \\
             & = & k(k-1) D(H,k;x,y)/k^2 \frac{ |Aut(H,k;x,y)|}{|Aut(H,k;xy)|} \\
             & = & k(k-1) D(H,k;x,y)/2k^2
\end{eqnarray*}
where the last equation follows from the fact that $|Aut(H;x,y)| =  |Aut(H; x\rightarrow y, y \rightarrow x)|$ because $Aut(H; x\rightarrow y, y \rightarrow x)$ is a coset of $Aut(H;x,y)$ and so  $|Aut(H;xy)| = 2 |Aut(H;x,y)|$.  Finally, since $D(H,k;xy) = D_1(H,k;xy) + D_2(H,k;xy)$ the formula for $D_1(H,k;xy)$ in the lemma follows. \qed

\begin{theorem}
\label{s-vertex-2}  
Let $v$ be an $s$-vertex in $T_G$ and $\{x,y\}$ be the separating pair associated with $v$.  If it exists, let $w_x$ denote the child of $v$ that is a $c$-vertex associated with $x$.  Similarly, if it exists, let $w_y$ denote the child of $v$ that is a $c$-vertex associated with $y$.  Suppose when all the graphs in $\mathcal{G} = \{ G(T_w):w \mbox{ is a $t$-vertex and a child of $v$ in $T_G$} \}$ are fixed at $x$ and $y$,  there are $g$ isomorphic classes and the $i$th isomorphic class has $m_i$ copies of the connected graph $G_i$; i.e. $\mathcal{G} = m_1 G_1 \cup m_2 G_2 \cup \hdots \cup m_g G_g$. If $Aut(G(T_v);xy) = Aut(G(T_v);x,y)$ then $D(G(T_v),k;xy) = D(G(T_v),k; x,y)$. Otherwise,  $Aut(G(T_v); x \rightarrow y, y \rightarrow x) \not = \emptyset$  so $G(T_{w_x}) \cong  G(T_{w_y})$, and 
 
$$D(G(T_v),k;xy) = [ D(G(T_v),k;x,y) - B(G(T_v),k;x,y)]/2$$ 
where $B(G(T_v),k;x,y)$ equals
 $$ k  \max\{D(G(T_{w_x}),k;x)/k,1\}  \prod_{i=1}^g  \sum_{l = 0}^{\lfloor m_i/2 \rfloor} \binom{D_1(G_i,k;xy)/k}{l} \binom{ [D(G_i,k;x,y) - 2D(G_i,k;xy)]/k   }{m_i - 2l}. $$ 
\end{theorem}

\noindent Proof:   When $Aut(G(T_v);xy) = Aut(G(T_v); x,y)$, $\mathcal{L}(G(T_v),k;xy) = \mathcal{L}(G(T_v),k;x,y)$ and so it follows that   $D(G(T_v),k;xy) = D(G(T_v),k; x, y)$.  Otherwise, $Aut(G(T_v); x \rightarrow y, y \rightarrow x) \not = \emptyset$.  Hence, $G(T_{w_x}) \cong G(T_{w_y})$ and  $Aut(G_i; xy) \not = Aut(G_i;x,y)$ for $i= 1, \hdots, g$.  In computing $D(G(T_v),k;x,y)$, we noted that there are three sets of  parameters that describe the equivalence classes of $\mathcal{L}(G(T_v),k;x,y)$: (i) the labels assigned to $x$ and $y$, (ii) the equivalence classes of $\mathcal{L}(G(T_{w_x}),k;x)$ that contain the labelings of $G(T_{w_x})$ and $G(T_{w_y})$, and (iii) the set of $m_i$ equivalence classes of $\mathcal{L}(G_i,k;x,y)$ that contain the labelings of the $m_i$ copies of $G_i$ for $i=1, \hdots, g$.  We shall extend them to characterize the equivalence classes in $\mathcal{B}(G(T_v),k;x,y)$ -- i.e., the equivalences classes of $\mathcal{L}(G(T_v),k;x,y)$ whose labelings are preserved by some automorphism in $Aut(G(T_v); x \rightarrow y, y \rightarrow x)$. 


\begin{claim}
 An equivalence class of $\mathcal{L}(G(T_v),k;x,y)$ belongs to $\mathcal{B}(G(T_v),k;x,y)$ if and only if 

\noindent (i) the labels assigned to $x$ and $y$ are the same for every labeling in the class, \\
 (ii) the equivalence classes that contain the labelings of $G(T_{w_x})$ and $G(T_{w_y})$ are the same, and  \\ 
(iii) for $i=1, \hdots, g$, the set of $m_i$ equivalence classes of $\mathcal{L}(G_i,k;x,y)$ that contain the labelings of the $m_i$ copies of $G_i$  can be partitioned into $l_i$ groups of size $2$ and $m_i - 2l_i$ groups of size $1$ for some $0 \leq l_i \leq \lfloor m_i/2 \rfloor$ so that the pairs of equivalence classes that belong to a group of size $2$ are partners in $\mathcal{L}(G_i,k;x,y)$ and the equivalence classes that belong to a group of size $1$ are in  $\mathcal{B}(G_i,k;x,y)$.
\end{claim}

\noindent Proof of claim:  Let $\phi \in \mathcal{L}(G(T_v),k;x,y)$ belong to an equivalence class that satisfies conditions (i), (ii) and (iii) above. For $i=1, \hdots, g$, denote the $m_i$ copies of $G_i$ as $G_{i,1}, \hdots, G_{i,m_i}$.  Without loss of generality,  assume that the equivalence classes that contain the labelings of $G_{i,2j-1}$ and $G_{i,2j}$ are partners in $\mathcal{L}(G_i,k;x,y)$  for $j=1, \hdots, l_i$ and the equivalence classes that contain the labelings  of each of the remaining copies of $G_i$ belong to $\mathcal{B}(G_i,k;x,y)$. Thus, condition (iii) implies that there are automorphisms~\footnote{Technically, we are referring to the automorphisms implicitly defined by the isomorphisms that maps $(G_{i,2j-1}, \phi)$  to $(G_{i,2j}, \phi)$ and vice versa since $G_{i,2j-1}$ and $G_{i,2j}$ are copies of $G_i$.  We shall keep this usage throughout the proof for ease of discussion.}  
 in $Aut(G_i; x\rightarrow y, y \rightarrow x)$ that map $(G_{i,2j-1}, \phi)$  to $(G_{i,2j}, \phi)$  and vice versa for $j=1, \hdots, l_i$; similarly, there is also some automorphism in $Aut(G_i; x\rightarrow y, y \rightarrow x)$ that preserves $(G_{i,j}, \phi)$ for $j = 2l_i + 1, \hdots, m_i$.  Furthermore, condition (ii) implies that there are some automorphisms in $Aut(G(T_{w_x});x)$ that map $(G(T_{w_x}), \phi)$ to $(G(T_{w_y}), \phi)$ and vice versa.  Combining these automorphisms, we conclude that some automorphism in $Aut(G(T_v); x\rightarrow y, y \rightarrow x)$ preserves $\phi$; that is, the equivalence class that contains $\phi$ belongs to $\mathcal{B}(G(T_v),k;x,y)$.

On the other hand, suppose $\phi \in \mathcal{L}(G(T_v),k;x,y)$ and some automorphism  $ \pi \in Aut(G(T_v); x\rightarrow y, y \rightarrow x)$ preserves $\phi$.  Clearly, $\pi$ maps $x$ to $y$ and vice versa, and so condition (i) is true. It also maps $G(T_{w_x})$ to $G(T_{w_y})$ and vice versa, and so condition (ii) is true.  For $i=1, \hdots, g$, if $\pi$ fixes $G_{i,j}$  then $\phi$ must have assigned $G_{i,j}$ a labeling that destroys all automorphisms in $Aut(G_i;x,y)$ since $\phi \in \mathcal{L}(G(T_v),k;x,y)$ but is still preserved by some automorphism in $Aut(G_i; x \rightarrow y, y \rightarrow x)$ since  $ \pi \in Aut(G(T_v); x\rightarrow y, y \rightarrow x)$.  That is, the equivalence class that contains the labeling of $G_{i,j}$ belongs to $\mathcal{B}(G_i,k;x,y)$.  If $\pi$ maps $G_{i,j}$ to $G_{i,j'}$, $j \not = j'$ then there is some automorphism in $Aut(G_i; x\rightarrow y, y \rightarrow x)$ that maps $(G_{i,j}, \phi)$ to $(G_{i,j'}, \phi)$ so that the equivalence classes that contain the labelings assigned by $\phi$ to  $G_{i,j}$ and $G_{i,j'}$ are partners in $\mathcal{L}(G_i,k;x,y)$.  Finally,  if   $\pi$ maps $G_{i,j}$ to $G_{i,j'}$ and $G_{i,j'}$ to $G_{i,j''}$ where $j, j'$ and $j''$ are distinct, then  the equivalence classes that contain the labelings assigned by $\phi$ to $G_{i,j}$ and $G_{i,j'}$ are partners  and those of $G_{i,j'}$ and $G_{i,j''}$ are partners as well.  But this implies that the equivalence classes that contain the labelings assigned by $\phi$ to $G_{i,j}$ and $G_{i,j''}$ are exactly the same, contradicting the assumption that $\phi$ assigned inequivalent labelings to the $m_i$ copies of $G_i$ since $\phi \in \mathcal{L}(G(T_v),k;x,y)$.  Thus, condition (iii) must be true. \qed
\smallskip

Combining conditions (i) and (iii), we note that each group of size $2$ in condition (iii) corresponds to an equivalence class of $\mathcal{L}_1(G_i,k;xy)$.  Moreover, because  no two  of the equivalence classes that contain the labelings of the $m_i$ copies of $G_i$  are identical,  distinct groups of size $2$ correspond to distinct equivalence classes of $\mathcal{L}_1(G_i,k;xy)$.    Using the claim, we can now compute $B(G(T_v),k;x,y)$. An equivalence class in $\mathcal{B}(G(T_v),k;x,y)$ is defined by (i) the label assigned to $x$ and $y$, (ii) the equivalence class of $\mathcal{L}(G(T_{w_x}),k;x)$ that  contain the labelings of $G(T_{w_x})$ and $G(T_{w_y})$, (iii) the set of $l_i$ equivalence classes of   $\mathcal{L}_1(G_i,k;x,y)$ and the set of $m_i-2l_i$ equivalence classes in $\mathcal{B}(G_i,k;x,y)$ which together contain the labelings of the $m_i$ copies of $G_i$ for $i=1, \hdots, g$.  
There are $k$ ways of assigning the same labels to $x$ and $y$, say $l$.  Once $l$ is fixed, there are $D(G(T_{w_x}),k;x)/k$ choices for the equivalence class of (ii), $\binom{D_1(G_i,k;x,y)/k}{l_i}$ choices for the set of $l_i$ equivalence classes and $\binom{B(G_i,k;x,y)/k}{m_i-2l_i}$ choices for the set of $m_i-2l_i$ equivalence classes of (iii).  Thus, $B(G(T_v),k;x,y)$ equals  



$$ k  \binom{D(G(T_{w_x}),k;x)/k}{1} \prod_{i=1}^g \sum_{l = 0}^{\lfloor m_i/2 \rfloor} \binom{D_1(G_i,k;xy)/k}{l} \binom{ [D(G_i,k;x,y) - 2D(G_i,k;xy)]/k   }{m_i - 2l}. $$

Once $B(G(T_v),k;x,y)$ has been computed, we can determine $D(G(T_v),k;xy)$ using Lemma~\ref{s-vertex-lemma-1}. \qed

\medskip 

An important implication of Theorems~\ref{s-vertex-1} and \ref{s-vertex-2} and Lemma~\ref{s-vertex-lemma-2}  is that both $D(G(T_v),k;x,y)$ and $D(G(T_v),k;xy)$ can be computed once the values of $D(G(T_{w_x}),k; x)$, $D(G(T_{w_y}),k; y)$, $D(G_i,k;x,y)$, and $D(G_i,k;xy)$ for $i=1, \hdots, g$ are known.  
\medskip

\begin{figure}
\label{t-vertex-example}
\centering
\epsfig{file=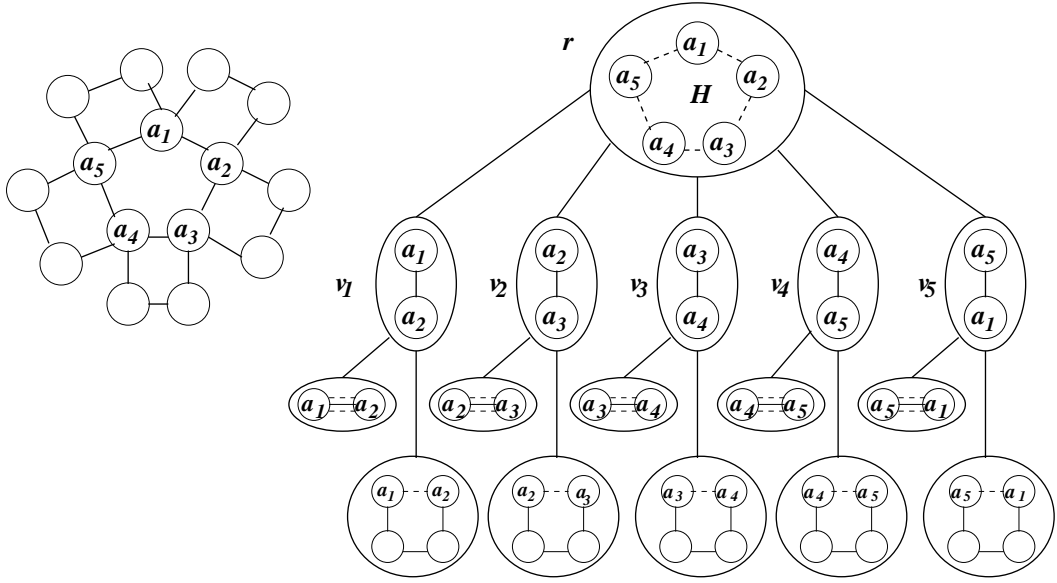, width=5.5in}
\caption{Graph $G$ and its tree decomposition.  Again, $r = r(T_G)$.   Note that $D(G(T_r),k; H) = D(G,k)$.} 
\end{figure}

Let us now consider the case when $v$ is a $t$-vertex.  Let $H$ be the triconnected component associated with $v$.  We need to solve for  $D(G(T_v),k; H,A)$ where $A$ is the structure associated with the parent of $v$ in $T_G$ (if the parent exists). Our goal is to create a formula for $D(G(T_v),k; H,A)$ that is dependent on $H$ and the values of $D(G(T_w),k;*)$ only, where $w$ a child of $v$ in $T_G$ and $*$ is the structure associated with $w$,  so that the formula can be computed efficiently.   Our approach follows Section~\ref{section-PIE} closely; the difference is that in our current setting $G(T_v)$ is made up of $H$ together with components hanging off of the cut vertices and separating pairs of $G$ in $H$ whereas in Section~\ref{section-PIE} we only dealt with the graph $H$. We demonstrate our approach by solving for $D(G(T_v),k;H)$; others can be solved similarly.  To aid us in  our discussion, we shall use the graph in Figure 6 for illustration.  In particular, we will compute for  $D(G(T_r),k; H)$, which equals $D(G,k)$.   


Let $\mathcal{C}_H$ contain the cut vertices of $G$ in $H$ whose corresponding vertices in $T_G$ are children of $v$.   Let $\mathcal{S}_H$ contain the  separating pairs of $G$ in  $H$ used in creating $T_G$ whose corresponding vertices in $T_G$ are children of $v$.   When  $a \in \mathcal{C}_H$  and $w$ is the child of $v$ that is associated with $a$, we shall refer to $G(T_w)$ as $G_a$ for ease of notation. We do the same for each pair $\{x,y\} \in \mathcal{S}_H$.    We begin by considering $\mathcal{L}^*(G(T_v),k)$, the set that contains all the $k$-labelings $\phi$ of $G(T_v)$ so that  $\phi$ when restricted to $G_a$  belongs to $\mathcal{L}(G_a,k;a)$ for every cut vertex $a \in \mathcal{C}_H$, and  $\phi$ when restricted to $G_{x,y}$ belongs to $\mathcal{L}(G_{x,y},k;x,y)$ for every separating pair $\{x,y\} \in \mathcal{S}_H$.    For instance, a labeling that assigns all the nodes of $G$ in Figure 6 the same label belongs to $\mathcal{L}^*(G(T_r),k)$ because every $G_{a_i, a_{i+1}}$ fixed at $a_i$ and $a_{i+1}$ has no non-trivial automorphisms. 
Clearly, every labeling in $\mathcal{L}(G(T_v),k;H)$ also belongs to $\mathcal{L}^*(G(T_v),k)$; otherwise, some nontrivial automorphism in $Aut(G(T_v);H)$ is not destroyed.  We shall use PIE to weed out the labelings in $\mathcal{L}^*(G(T_v),k)$ that are not part of $\mathcal{L}(G(T_v),k;H)$.  

Let $Aut_v(H)$ consist of all the automorphisms in $Aut(G(T_v);H)$ when restricted to $H$.  Suppose $P \subseteq Aut_v(H)$.  Let $N_{\ge}(P)$ denote the number of labelings $\phi$ in $\mathcal{L}^*(G(T_v),k)$ so that, for every $\sigma \in P$, some extension of $\sigma$ in $Aut(G(T_v);H)$ preserves $\phi$.  Define $N_=(P)$ similarly except that aside from the automorphisms in $P$ no other automorphism in $Aut_v(H)$ has extensions that preserve $\phi$.  We now state the formulas for $L(G(T_v),k;H)$ and $D(G(T_v),k;H)$.

\begin{theorem}
Let $v$ be a $t$-vertex in $T_G$ and  $H$ be the triconnected component associated with $v$.
 Let $Aut_v(H)$ consist of the automorphisms in $Aut(G(T_v);H)$ when restricted to $H$, and $\sigma_0$ be the identity automorphism in $Aut_v(H)$.  Then
$$L(G(T_v),k;H) = N_=(\{ \sigma_0 \}) =  \sum_{ \{ \sigma_0 \} \subseteq P \subseteq Aut_v(H)} (-1)^{|P|-1} N_\ge(P) $$
and 
$$ D(G(T_v),k;H) = \frac{1}{|Aut_v(H)|} \sum_{ \{ \sigma_0 \} \subseteq P \subseteq Aut_v(H)} \frac{(-1)^{|P|-1} N_\ge(P)}{\prod_{a \in \mathcal{C}_H}|Aut(G_a;a)| \prod_{\{x,y\}\in \mathcal{S}_H} |Aut(G_{x,y};x,y)|}. $$
\label{t-vertex-1}
\end{theorem}

\noindent Proof:  It is straightforward to verify that $ N_=(\{ \sigma_0 \})  =  \sum_{ \{ \sigma_0 \} \subseteq P \subseteq Aut_v(H)} (-1)^{|P|-1} N_\ge(P)$ is the number of labelings in  $\mathcal{L}^*(G(T_v),k)$ that is preserved by  some extension of $\sigma_0$ and by no other automorphism in $Aut_v(H)$.   But by the way we defined $\mathcal{L}^*(G(T_v),k)$,  if some extension of $\sigma_0$  preserves a labeling of $\mathcal{L}^*(G(T_v),k)$, that extension must be the identity automorphism of $G(T_v)$.  The first equation of the theorem follows.  Now, $|Aut(G(T_v);H)| = |Aut_v(H)| \times \prod_{a\in \mathcal{C}_H} |Aut(G_a;a)| \times  \prod_{\{x,y\}\in \mathcal{S}_H} |Aut(G_{x,y};x,y)|$ since every automorphism in $Aut(G(T_v); H)$ can be decomposed into an automorphism in $Aut_v(H)$ and automorphisms of the connected components that hang off of $H$ fixed at the cut vertices or separating pairs that connect them to $H$.  Hence, dividing $ L(G(T_v),k;H)$ by $|Aut(G(T_v);H)|$ gives us the second equation of the theorem. \qed

\begin{corollary}
\label{t-vertex-corollary}
When $Aut_v(H)$ is cyclic, dihedral, or isomorphic to $Z_t \times Z_2$ or $D_t \times Z_2$ for some integer $t$, the formulas in Theorems~\ref{cyclic}, \ref{dihedral} and \ref{crossproduct} for $N_=(\{\sigma_0 \})$ where $\sigma_0$ is the identity automorphism still apply.
\end{corollary}

Suppose, instead of $D(G(T_v),k;H)$, we wish to compute  $D(G(T_v),k; H, a)$ (or $D(G(T_v),k; H, x,y)$ or $D(G(T_v),k; H, xy)$).  Let $\Gamma$ denote the subgroup of $Aut_v(H)$ that fixes $a$ (or $x$ and $y$, or $xy$). By replacing $Aut_v(H)$ with $\Gamma$,  the formula in Theorem~\ref{t-vertex-1} still holds.

Next, we describe a method for computing $N_\ge(P)$, $P \subseteq Aut_v(H)$.  For each $\sigma \in P$,  let $H_{\sigma}$ be the graph whose vertex set is $V(H)$ and edge set is $\{(v, \sigma(v)): v \in V(H)\}$. Let $SP_{\sigma}$ be the graph whose vertex set consists of $xy$ and $yx$ whenever $\{x,y\} \in \mathcal{S}_H$ and whose edge set is $\{(xy, \sigma(x)\sigma(y)), (yx, \sigma(y)\sigma(x)): \{x,y\} \in \mathcal{S}_H\}$.  In $\cup_{\sigma \in P} SP_{\sigma}$,  let $SP_P(xy)$ denote the component  that contains the vertex $xy$.   
Notice that  if $a$ and $a'$ are part of the same component in $\cup_{\sigma \in P} H_{\sigma}$, then both $a$ and $a'$ are cut vertices or both are not; and when both of them are, then $G_a \cong G_{a'}$.  Similarly, when $xy$ and $x'y'$ are part of the same component in $\cup_{\sigma \in P} SP_{\sigma}$, $G_{x,y} \cong G_{x',y'}$.  Finally,  when $yx \in SP_P(xy)$ then for every $x'y' \in SP_P(xy)$, $y'x' \in SP_P(xy)$ as well.  The following can easily be verified: $\phi \in \mathcal{L}^*(G(T_v),k)$ is counted  in $N_\ge(P)$ if and only if  (i) whenever two vertices are part of the same component in $\cup_{\sigma \in P} H_{\sigma}$, $\phi$ assigns them the same label; (ii) whenever two cut vertices $a$ and $a'$ are part of the same component in $\cup_{\sigma \in P} H_{\sigma}$, $\phi$ when restricted to $G_a$ and $G_{a'}$ belong to the same equivalence class of $\mathcal{L}(G_a,k;a)$; and (iii) when $xy$ and $x'y'$ are part of the same component in $\cup_{\sigma \in P} SP_{\sigma}$, $\phi$ when restricted to $G_{x,y}$ and $G_{x',y'}$ belong to the same equivalence class of $\mathcal{L}(G_{x,y},k;x,y)$. 

In $SP_{\sigma}$, we have chosen to represent the separating pair $\{x,y\}$ as two vertices $xy$ and $yx$ to capture situations in which an automorphism in $Aut_v(H)$ maps $x$ to $y$ and $y$ to $x$.   However, such a representation can introduce redundancies in the sense that two different components in $SP_\sigma$ may be capturing the same relationships between the same sets of separating pairs.  As such, we   shall say that a collection of components  $\{SP_P(x_iy_i), i = 1, \hdots g\}$ forms a {\it partition} of $\cup_{\sigma \in P} SP_{\sigma}$ if for every separating pair $\{x,y\} \in \mathcal{S}_H$ exactly one  component in the collection contains $xy$ or $yx$ or both.  In Figure 6, let $\pi_{\rho} \in Aut_r(H)$ be the rotation that maps $a_i$ to $a_{i+1}$, $i = 1, \hdots, 5$,  and $\pi_{\tau} \in Aut_r(H)$ be the reflection that fixes $a_1$.  A partition for $SP_{\pi_{\rho}}$ contains only one component while a partition for $SP_{\pi_{\tau} }$ contains three components -- e.g., $SP(a_1a_5), SP(a_5a_4), SP(a_4a_3)$ where  $a_3a_4 \in SP(a_4a_3)$.    


\begin{theorem}
\label{t-vertex-2}
Let $\{\sigma_0 \} \subseteq P \subseteq Aut_v(H)$.  Suppose $\cup_{\sigma \in P} H_{\sigma}$ has $t$ components where the $i$th component contains the vertex $a_i$,  and the collection $\{SP_P(x_iy_i), i = 1, \hdots g\}$ forms a {\it partition} of $\cup_{\sigma \in P} SP_{\sigma}$. Let $\kappa(a_i) = D(G_{a_i},k;a_i)$ if $a_i$ is a cut vertex and is equal to $k$ otherwise.  Let $\kappa(x_i,y_i) =  [D(G_{x_i,y_i},k; x_i,y_i) - 2 D(G_{x_i,y_i},k; x_iy_i)]/k$ if  $y_ix_i \in SP_P(x_iy_i)$ and is equal to $D(G_{x_i,y_i},k;x_i,y_i)/k^2$ otherwise. Then  

$$ N_\ge(P) = \prod_{i=1}^t \kappa(a_i) \times  \prod_{i=1}^g \kappa(x_i,y_i) \times \prod_{a \in \mathcal{C}_H}|Aut(G_a;a)| \times \prod_{\{x,y\}\in \mathcal{S}_H} |Aut(G_{x,y};x,y)|.$$

\end{theorem}

\noindent Proof: To create a labeling $\phi$ that is counted  in $N_\ge(P)$, we do the following. (1) If $a_i$ is not a cut vertex, pick a label that will be assigned to it and all the vertices in the same component as $a_i$ in $\cup_{\sigma \in P} H_{\sigma}$. (2) If $a_i$ is a cut vertex, pick an equivalence class of $\mathcal{L}(G_{a_i},k; a_i)$ that will contain the labelings of $G_{a_i}$ and all the $G_u$'s where $u$ and $a_i$ are in the same component of $\cup_{\sigma \in P} H_{\sigma}$.  Then for $a_i$ and each of the vertices $u$, pick a labeling from the equivalence class just chosen. (3) For each $\{x_i,y_i\}$, pick an equivalence class of $\mathcal{L}(G_{x_i,y_i},k; x_i,y_i)$ that will contain the labelings of $G_{x_i,y_i}$ and all $G_{u,w}$'s where $uw$ and $x_iy_i$ are in the same component of $\cup_{\sigma \in P} SP_{\sigma}$.  Additionally, this equivalence class must respect the labels that have already been assigned to $x_i$ and $y_i$ in step (2).  Then for $\{x_i,y_i\}$ and each $\{u,w\}$, pick a labeling from the equivalence class just chosen.

There are $k$ ways of doing step (1) and $D(G_{a_i},k;a_i) \times |Aut(G_{a_i};a_i)|^{j_i}$ ways of doing step (2) where $j_i$ is the number of vertices in the same component as $a_i$.  Thus, there are $\prod_{i=1}^t \kappa(a_i)  \prod_{a \in \mathcal{C}_H} |Aut(G_a;a)|$ ways of doing steps (1) and (2) since whenever the $a_i$ and $u$ are in the same component of $\cup_{\sigma \in P} H_{\sigma}$, both vertices are cut vertices or both are not and $Aut(G_{a_i};a_i) \cong Aut(G_u;u)$.   

To do step (3), we need to differentiate between the case when $y_ix_i$ also belongs to $SP_P(x_iy_i)$ and when it doesn't.  In the former case, $\phi$ when restricted to $G_{x_i,y_i}$ must destroy all the automorphisms of $G_{x_i,y_i}$ when $x_i$ and $y_i$ are fixed but is preserved by some automorphism of the graph that maps $x_i$ to $y_i$ and vice versa.  In other words, the equivalence class containing the labeling belongs to $\mathcal{B}(G_{x_i,y_i},k;x_i,y_i)$.  Since the labels of $x$ and $y$ have already been chosen in steps (1) or (2)  (note that they had to be the same), from Lemma~\ref{s-vertex-lemma-1} there are exactly $[D(G_{x_i,y_i},k ; x_i, y_i) - 2D(G_{x_i,y_i}, k; x_iy_i)]/k$ equivalence classes to choose from in step (3).  On the other hand, when $y_ix_i$ does not belong to $SP_P(x_iy_i)$, $\phi$ when restricted to $G_{x_i,y_i}$ must simply belong to $\mathcal{L}(G_{x_i,y_i},k;x_i,y_i)$ and so once the labels of $x_i$ and $y_i$ have been chosen in steps (1) and (2), there are exactly $D(G_{x_i,y_i},k;x_i,y_i)/k^2$ equivalence classes to choose from in step (3).  Finally, when the equivalence classes have been selected, then there are $|Aut(G_{x_i,y_i};x_i,y_i)|^{j_i}$ labelings that can be assigned to the $G_{u,w}$'s, $uw \in SP_P(x_iy_i)$ where $j_i$ is the number of distinct separating pairs in $SP(x_iy_i)$.  There are $\prod_{i=1}^g \kappa(x_i,y_i) \prod_{\{x,y\}\in \mathcal{S}_H} |Aut(G_{x,y};x,y)|$ ways of doing step (3) because again whenever $x_iy_i$ and $uw$ are in the same component in $\cup_{\sigma \in P} SP_{\sigma}$, $Aut(G_{x_i,y_i}; x_i,y_i) \cong Aut(G_{u,w}; u,w)$.  The theorem follows. \qed
\medskip

\medskip


Using the formula for $N_{\ge}(P)$ above, we can now simplify the second formula in  Theorem~\ref{t-vertex-1} as  
$$ D(G(T_v),k;H) = \frac{1}{|Aut_v(H)|} \sum_{ \{ \sigma_0 \} \subseteq P \subseteq Aut_v(H)} (-1)^{|P|-1} \prod_{i=1}^t \kappa(a_i) \prod_{i=1}^g \kappa(x_i,y_i). $$

Hence, $D(G(T_v),k;H)$ (and $D(G(T_v),k;H,a)$,  $D(G(T_v),k; H, x,y)$, $D(G(T_v),k; H, xy)$) can be computed once  the automorphisms in $Aut_v(H)$ and the values of $D(G(T_w),k;*)$ are known, where $w$ a child of $v$ in $T_G$ and $*$ is the structure associated with $w$. 


\smallskip

\noindent {\it FindDist($G,k$).}  Let us now describe our algorithm {\it FindDist($G,k$)} for computing $D(G,k)$.   First, construct $T_G$ and root it at $r(T_G)$.  Then, for $i = \mbox{depth$(T_G)$ to $0$}$,   do the following for each vertex $v$ at depth $i$.  When $v$ is a $c$-vertex, compute $D(G(T_v),k;a)$ where $a$ is the cut vertex associated with $v$ using Theorem~\ref{c-vertex}. When $v$ is an $s$-vertex, compute $D(G(T_v),k;x,y)$ and $D(G(T_v),k;xy)$ where $\{x,y\}$ is the separating pair associated with $v$ using Theorems~\ref{s-vertex-1} and \ref{s-vertex-2}. When $v$ is a $t$-vertex and $H$ is the triconnected component associated with $v$, compute $D(G(T_v),k; H,*)$ where $*$ is the structure associated with $p(v)$, the parent of $v$, if it exists.  (That is, if $p(v)$ is a $c$-vertex,  compute $D(G(T_v),k; H, a)$ where $a$ is the cut vertex associated with $p(v)$; if $p(v)$ is an $s$-vertex, compute $D(G(T_v),k; H, x,y)$ and $D(G(T_v),k; H, xy)$ where $\{x,y\}$ is the separating pair associated with $p(v)$; if $v$ has no parent, compute $D(G(T_v),k;H)$.)   Do all  computations  according to Theorem~\ref{t-vertex-1}, Corollary~\ref{t-vertex-corollary} and Theorem~\ref{t-vertex-2}.  Finally, if the root node $r(T_G)$ is a $c$-vertex or a $t$-vertex, return the value computed at depth $0$; otherwise, $r(T_G)$ is an $s$-vertex, return $D(G(T_v),k;xy)$. 

It is not obvious that at every iteration the algorithm can compute $D(G(T_v),k;*)$  based on the values obtained for $v$'s children in the previous iterations. We shall now show that this, in fact, is the case.  


\noindent $\bullet$ {\it When $v$ is a $c$-vertex.}  Let $a$ be the cut vertex associated with $v$.  From Theorem~\ref{c-vertex}, $D(G(T_w),k;a)$ is needed.  If $w$ is an $s$-vertex and associated with some separating pair $\{a,b\}$, $D(G(T_w),k;a,b)$ was computed in the previous iterations. But  from Lemma~\ref{autc-vertex}, $Aut(G(T_w);a,b) = Aut(G(T_w);a)$; i.e., a labeling of $G(T_w)$ destroys all automorphisms in $Aut(G(T_w);a,b)$ if and only if it destroys all automorphisms in $Aut(G(T_w);a)$.  Hence, $D(G(T_w), k;a,b) = D(G(T_w), k;a)$. If $w$ is a $t$-vertex, $w$ is associated with some triconnected component $H$, and so $D(G(T_w),k;H,a)$ was computed in the previous iteration.  From Lemma~\ref{autc-vertex}, $Aut(G(T_w);H,a) = Aut(G(T_w);a)$ and,  consequently, $D(G(T_w), k; H,a) = D(G(T_w), k;a)$.

\noindent $\bullet$ {\it When $v$ is an $s$-vertex.}  Let $\{x,y\}$ be the separating pair associated with $v$.  This time, according to Theorems~\ref{s-vertex-1} and \ref{s-vertex-2}, if $w$ is a $c$-vertex associated with $x$ (or $y$), $D(G(T_w),k;x)$ (or $D(G(T_w),k;y)$) is needed, and was clearly computed in the previous iterations by the algorithm.  On the other hand, if $w$ is a $t$-vertex whose associated triconnected component is $H$,  $D(G(T_w),k;x,y)$ and $D(G(T_w),k;xy)$ are needed.  Now, $D(G(T_w),k; H,x,y)$ and $D(G(T_w),k; H,xy)$ were computed in the previous iterations.  But from Lemma~\ref{auts-vertex}, $Aut(G(T_w); H,x,y) = Aut(G(T_w); x,y)$ and $Aut(G(T_w); H,xy) = Aut(G(T_w); xy)$. It follows that $D(G(T_w),k; H,x,y) = D(G(T_w),k; x,y)$ and $D(G(T_w),k; H,xy) = D(G(T_w),k; xy)$.  

\noindent $\bullet$ {\it When $v$ is a $t$-vertex.}  Let $H$ be the triconnected component associated with $v$.  According to Theorem~\ref{t-vertex-2}, if $w$ is a $c$-vertex associated with $a$, $D(G(T_w),k;a)$ is needed and if $w$ is an $s$-vertex associated with the separating pair $\{x,y\}$, $D(G(T_w),k;x,y)$ and $D(G(T_w),k;xy)$ are needed.   All these values were computed in the previous iterations.
\smallskip

From Theorems~\ref{c-vertex}, \ref{s-vertex-1}, \ref{s-vertex-2}, \ref{t-vertex-1}, \ref{t-vertex-2}, we know that all the $D(G(T_v),k;*)$ values computed by the algorithm are correct.  Now, the algorithm returned the value $D(G(T_v),k; A)$ where $v = r(T_G)$ and $A$ is the structure associated with $r(T_G)$.   But $G(T_{r(T_G)}) = G$, and, according to Lemma~\ref{autroot}, $Aut(G;A) = Aut(G)$.  It follows that the algorithm returned $D(G,k)$.  


\begin{theorem}
When $G$ is a connected graph, {\it FindDist($G,k$)} returns the value $D(G,k)$.
\end{theorem}

\smallskip

\noindent {\it Example 1.}  Consider the graph and its tree decomposition  in Figures 2 and 3.  Using the formulas given in this section, it is easy to verify the following:  $D(G(T_{v_3}),k;e) = k^5$, $D(G(T_{v_1}),k;e,j) = k^3(k-1)/2$ and $D(G(T_{v_1}),k;ej) = k^2(k-1)^2/4$.   According to Theorem~\ref{s-vertex-1},
\begin{eqnarray*}
D(G(T_r),k;e,j) & = & k^2 [D(G(T_{v_3}),k; e)/k]^2 \binom{D(G(T_{v_1}),k;e,j)/k^2}{2} \\ 
   & = & k^2 (k^4)^2 \binom{k(k-1)/2}{2} \\
   & = & (k+1) k^{11} (k-1)(k-2)/8. 
\end{eqnarray*}

Next, let us compute $B(G(T_r),k;e,j)$.  We have $B(G(T_{v_1}),k;e,j) = D(G(T_{v_1}),k;e,j) - 2D(G(T_{v_1}),k;ej) = k^2(k-1)/2$.
 Notice that $D_1(G(T_{v_1}),k;ej) = 0$ because any labeling that assigns $e$ and $j$ the same label cannot destroy the automorphism that maps $e$ to $j$, $j$ to $e$, $f$ to itself, and $g$ to itself.  According to Theorem~\ref{s-vertex-2}
\begin{eqnarray*}
B(G(T_r),k;e,j) & = & k \:D(G(T_{v_3}),k; e)/k \: \sum_{l=0}^1 \binom{D_1(G(T_{v_1}),k;e,j)/k}{l} \binom{B(G(T_{v_1}),k; e,j)/k}{2-2l} \\
   & = & k \: k^4 \: \sum_{l=0}^1 \binom{0}{l}\binom{k(k-1)/2}{2-2l} \\
   & = & k^5 \left[ \binom{0}{0}\binom{k(k-1)/2}{2} + \binom{0}{1} \binom{k(k-1)/2}{0} \right] \\
   & = & (k+1)k^6 (k-1)(k-2)/8. 
\end{eqnarray*}

Consequently, $D(G,k) = D(G(T_r),k;ej) = [D(G(T_r),k;e,j) - B(G(T_r),k;e,j)]/2 = (k+1)k^6(k-1)(k-2)(k^5-1)/16 .$  Since $D(G,1) = D(G,2) = 0$ and $D(G,3) > 0$,  $D(G) = 3$. 
\medskip

\noindent {\it Example 2.}  This time, consider the graph in Figure 6 and let us determine $D(G(T_r),k;H)$. Since $Aut_r(H) \cong D_5$, we can make use of the computations we made in the example after Theorem~\ref{dihedral}. 
It is easy to verify that $D(G_{a_i,a_{i+1}},k; a_i, a_{i+1}) = k^4$ and $D(G_{a_i,a_{i+1}},k; a_i a_{i+1}) = (k^4 - k^2)/2$ for $i=1, \hdots, 5$. Let $\alpha = \prod_{i=1}^5 |Aut(G_{a_i,a_{i+1}}; a_i, a_{i+1})|$. 
From Theorem~\ref{t-vertex-2}, we have $N_\ge(\emptyset)/\alpha = k^5 \times (k^2)^5 = k^{15}$, $N_\ge(\{\pi_{\rho}\})/\alpha = k \times k^2 = k^3$, $N_\ge(  \{ \pi_{\tau \rho^i} \})/\alpha = k^3 \times (k^2)^2 (k^4 - 2(k^4-k^2)/2)/k = k^8 $ and  $N_\ge(\{ \pi_{\tau \rho^i}, \pi_{\rho} \})/\alpha = k \times (k^4 - 2(k^4-k^2)/2)/k = k^2$.   Thus,
\begin{eqnarray*}
D(G(T_r),k; H) & = & \frac{1}{10} \left[ N_\ge(\emptyset)/\alpha - N_\ge(\{\pi_{\rho}\})/\alpha - \sum_{i=0}^{4}\left( N_\ge(  \{ \pi_{\tau \rho^i} \})/\alpha - N_\ge(\{ \pi_{\tau \rho^i}, \pi_{\rho} \})/\alpha \right)  \right] \\
               & = & \frac{1}{10} \left[k^{15} - k^3 - \sum_{i=0}^4 (k^8 - k^2) \right] \\
               & = & \frac{1}{10} \left[ k^{15} - 5 k^8 - k^3 + 5k^2 \right].
\end{eqnarray*}
Since $D(G,1) = 0$ and $D(G,2) > 0$, we conclude that $D(G) = 2$.
\medskip

In order for {\it FindDist($(G,k)$)} to run efficiently for a family of graphs, we note that a few ingredients are necessary.  It  must have an efficient graph isomorphism testing algorithm  to determine the isomorphism classes of $\mathcal{G}$ in Theorems~\ref{c-vertex}, \ref{s-vertex-1}, and \ref{s-vertex-2}.  There must also be an efficient algorithm that can determine the automorphisms of its triconnected components  which are  needed in Theorem~\ref{t-vertex-1}.  Finally, there must be a way to apply the PIE formula in Theorem~\ref{t-vertex-1} to its triconnected components in an efficient manner.  Since the family of planar graphs satisfy these criteria, we  can now proceed to prove the main result of the paper.  We first show that $D(G(T_v),k;*)$, when $v$ is a $t$-vertex, can be computed efficiently when the appropriate values are known.


\begin{lemma}
\label{t-vertex-planar}
Let $G$ be a connected planar graph on $n$ vertices and $k$ be a positive integer.  Let $v$ be a $t$-vertex in $T_G$, $H$ be the triconnected component on $n_H$ vertices associated with $v$,  and $\mathcal{S}_H$ contain the separating pairs of $G$ in $H$ used in the construction of $T_G$.  Suppose all the automorphisms in $Aut_v(H)$ and the values of $D(G(T_w),k;*)$, $w$ a child of $v$ in $T_G$,  are known.  Then $D(G(T_v),k;*)$ can be computed in $O(n_H^2 n^2 \log^2k (n_H + |\mathcal{S}_H|))$ time.
\end{lemma}

\noindent Proof:  As in Theorems~\ref{t-vertex-1} and \ref{t-vertex-2}, we will prove the theorem for $D(G(T_v),k;H)$; others can be argued similarly.  Let $\alpha = \prod_{a \in \mathcal{C}_H}|Aut(G_a;a)| \times \prod_{\{x,y\}\in \mathcal{S}_H} |Aut(G_{x,y};x,y)|$.
First, let us consider the time it takes to compute $N_{\ge}(P) / \alpha$.  Constructing $\cup_{\sigma \in P} H_{\sigma}$ and finding its $t$ connected components can be done in $O(|P| \times n_H)$.  Similarly, constructing $\cup_{\sigma \in P} SP_{\sigma}$ and finding a collection $\{SP_P(x_iy_i), i = 1, \hdots, g\}$ that forms a partition of $\cup_{\sigma \in P} SP_{\sigma}$ takes  $O(|P| \times |\mathcal{S}_H|)$ time.  Finally, since $N_{\ge}(P) \leq k^n$, every multiplication in $\prod_{i=1}^t \kappa(a_i) \times \prod_{i=1}^g \kappa(x_i y_i)$ takes at most $O(n^2 \log ^2k)$ time; that is, finding the said product takes $O((t+g) n^2 \log^2k) = O((n_H + |\mathcal{S}_H|)n^2 \log^2k)$ time.  Therefore, computing $N_{\ge}(P) / \alpha$  takes $O((n_H + |\mathcal{S}_H|)(|P| + n^2 \log^2k))$ time.

Since $G$ is a connected planar graph, each of its blocks $B$ are as well.  It is easy to verify that when the split operation is applied to $B$, the resulting split graphs remain planar.  Hence, the triconnected component associated with $v$ is  either a bond, a cycle or a triconnected planar graph, and $Aut_v(H)$ belongs to one of the groups mentioned in Fact~\ref{fact-lemma}.  From Corollary~\ref{t-vertex-corollary}, we noted we can apply the formulas in Theorems~\ref{cyclic}, \ref{dihedral} and \ref{crossproduct} for $N_=(\{\sigma_0 \})$.  Consequently, we can use  {\it TriconnectCount}$(H,k)$ to compute $D(G(T_v),k;H)$ by doing two  modifications -- we replace every occurrence of  $N_{\ge}(P)$ with  $N_\ge(P)/\alpha$ (keeping in mind that the $k$-labelings in $N_{\ge}(P)$  must all belong to $\mathcal{L}^*(G(T_v),k)$), and return the value $L(H,k)/|Aut_v(H)|$ instead of $L(H,k)$.  It is easy to check that the value returned is $D(G(T_v),k;H)$.  Applying the  same analysis in Theorem~\ref{triconnected-theorem}, the bottleneck of the algorithm is in computing the $O(n_H^2)$  $N_\ge(P)/ \alpha$ terms. Using the result from the previous paragraph and the fact that $|P| \leq \log n_H$, the runtime of the algorithm is $O(n_H^2(n_H + |\mathcal{S}_H|)(\log n_H + n^2 \log^2k))$ time. \qed
\medskip

Finding $D(G(T_v),k;*)$ can take significantly more time than finding $D(H,k)$ because  $N_\ge(P)$ in the former case has to be explicitly computed from the $D(G(T_w),k:*)$ values  whereas, in the latter case, $N_\ge(P)$ is always some power of $k$ which we precomputed ahead of time.

\begin{theorem}
Let $G$ be an $n$-vertex connected planar graph and $k$ be a positive integer.  {\it FindDist($G,k$)} can be implemented in $O(n^5 \log^2k)$ time.  Consequently, computing $D(G)$ takes $O(n^5 \log^3n)$ time.
\end{theorem}

\noindent Proof: To implement {\it FindDist($G,k$)}, we need to do  some preprocessing steps.  

\noindent 1. {\it Isomorphism testing.} For each $c$- and $s$-vertex, determine the isomorphism classes of the graphs in $\mathcal{G}$ fixed at the appropriate vertices as described in Theorems~\ref{c-vertex}, \ref{s-vertex-1}, and \ref{s-vertex-2}.  Additionally, for each $s$-vertex whose associated separating pair is $\{x,y\}$, determine if $G(T_{w_x})$ fixed at $x$ is isomorphic to $G(T_{w_y})$ fixed at $y$.  For each $t$-vertex, group together its children based on their type -- i.e., if they are $c$-vertices or $s$-vertices -- and then for each group determine the isomorphism classes of $\{G(T_w): \mbox{$w$ is a child of $v$ and a $c$-vertex (or an $s$-vertex)}\}$ fixed at the appropriate vertices.  Using the linear-time planar graph isomorphism testing algorithm, the isomorphism testing tasks at each $v$ can be done in $O(deg(v)^2 n)$ time where $deg(v)$ is the degree of vertex $v$ in $T_G$.  The total time to do all the isomorphism testing then takes $O(n^3)$ time since the size of $T_G$ is $O(n)$.    

\noindent 2. {\it Finding automorphisms.} For each $t$-vertex $v$ whose associated triconnected component is $H$, find all the automorphisms in $Aut(H)$ using the algorithm described in Section 3.1.    Apply the information obtained from the first preprocessing step to determine which of the automorphisms also belong to $Aut_v(H)$ (or $Aut_v(H; *)$ where $*$ depends on $p(v)$).  That is, for each $\pi \in Aut(H)$, verify that $\pi$ maps cut vertices and separating pairs to cut vertices and separating pairs respectively; moreover, $\pi$ also maps the subgraphs hanging from  these vertices and pairs to isomorphic structures.  Finding $Aut(H)$ takes $O(n_H^2)$ time where $n_H$ is the number of vertices in $H$.   Determining if $\pi \in Aut(H)$ also belongs to $Aut_v(H)$ takes $O(n_H + |\mathcal{S}_H|)$ time  so finding $Aut_v(H)$ takes  $O(|Aut(H)|(n_H + |\mathcal{S}_H|)) = O(n_H^2 + n_H |\mathcal{S}_H|)$ time.  Finding $Aut_v(H)$ for all $t$-vertices $v$ then takes $O(n^2)$ time since from Lemma~\ref{t-vertices-lemma} $\sum_H n_H = O(n)$ and $\sum_H |\mathcal{S}_H| = O(n)$.

Let us now consider {\it FindDist($G,k$)}.  Constructing $T_G$ takes $O(n^2)$ time. The number of arithmetic operations for computing $D(G(T_v),k;*)$ when $v$ is a $c$-vertex and an $s$-vertex is $O(deg(v))$ and $O(deg^2(v))$ respectively.   And since the values involved in each operation is at most $k^n$, the total amount of work at $v$ takes $O(deg^2(v) n^2 \log^2k)$ time.  Summing this up for all $c$- and $s$-vertices in $T_G$, the total work takes $O(n^4 \log^2k)$ time because the size of $T_G$ is $O(n)$.  When $v$ is a $t$-vertex,  according to Lemma~\ref{t-vertex-planar}, computing $D(G(T_v),k;*)$ takes $O(n_H^2 n^2 \log^2k (n_H + |\mathcal{S}_H|))$ time.  Thus, for all $t$-vertices, the total work is $O(n^5 \log^2k)$ time.  From our analysis, the bottleneck of {\it FindDist($G,k$)} is in processing the $t$-vertices.  The runtime in the theorem follows.  \qed

\begin{corollary}
Let $G$ be an $n$-vertex planar graph.  Then $D(G)$ can be computed in $O(n^5 \log^3n)$ time.
\end{corollary}

 \noindent Proof:  Using the linear time planar graph isomorphism testing  algorithm,  group together isomorphic connected components of $G$.  So suppose $G = m_1 G_1 \cup m_2 G_2 \cup \hdots \cup m_g G_g$ where each $G_i$ has $n_i$ vertices.  For each $G_i$, find $k_i$ so that $k_i = \min \{k: D(G_i,k) \ge m_i\}$.  This can be done by implementing  {\it FindDist$(G_i,k)$} $O(\log n)$ times.  Finally, $D(G) = \max_i \{k_i\}$ so $D(G)$ can be found in $O(\sum_{i=1}^g n_i^5 \log^2 n_i \log n) = O(n^5 \log^3 n)$ time. \qed

\section{Final Comments}


In this paper, we considered $D(G,k)$, the number of inequivalent distinguishing $k$-labelings of graph $G$.  We have applied the principle of inclusion/exclusion and developed recursive formulas to compute its value for a fixed $k$.  When a graph $G$ is planar, we showed that these techniques led to an algorithm for computing $D(G)$ that runs in time polynomial in the size of $G$.  There are other interesting aspects about $D(G,k)$ as well as noted in the next theorem.

\begin{theorem}
Let $G$ be a graph on $n$ vertices.  Then $D(G,k)$ is a polynomial in $k$ whose degree is $n$ and whose constant term is $0$.  If $G$ has no non-trivial automorphisms then $D(G,k) = k^n$; otherwise, the sum of the coefficients of $D(G,k)$ is $0$.   
\end{theorem}

\noindent Proof:  In formula~\ref{PIE},  $N_{\ge}(P)$  equals $k^n $ when 
$P = \{\pi_0\}$ and $k^t$, $t < n$, otherwise,  so $D(G,k)$ must be a polynomial in $k$ whose degree is $n$. Furthermore, $D(G,0) = 0$, so the constant term in $D(G,k)$ must be $0$.  When $G$ has no non-trivial automorphism, every $k$-labeling of $G$ is distinguishing and no two are equivalent; hence, $D(G,k) = k^n$.  When $G$ has some non-trivial automorphisms, $D(G,1) = 0$ so it must be the case that the sum of the coefficients of $D(G,k)$ is $0$.  \qed
\smallskip

We now call $D(G,k)$ the {\it distinguishing polynomial of $G$}.  An interesting research direction would be to study this polynomial along the lines of the more famous chromatic polynomial of $G$.
\smallskip

Next, consider Lemma~\ref{subgroup} and its implications.  According to the lemma,  if $S$ is the largest subset of $Aut(G)$ that preserves a $k$-labeling $\phi$ of $G$ then $S$ must be a subgroup of $Aut(G)$. This suggests that instead of considering all the subsets of $Aut(G)$ as we do in PIE, we should just consider the subgroup lattice $(\mathcal{S}, \leq)$ of $Aut(G)$ where $\mathcal{S}$ is the set that contains all the subgroups of $Aut(G)$. For each $S \in \mathcal{S}$, define $N_\ge(S)$ and $N_=(S)$ as we did in the PIE formulation. Since $N_\ge(S) = \sum_{S' \ge S} N_=(S')$ for every $S \in \mathcal{S}$, according to the principle of M\"obius inversion, the following must be true.

\begin{theorem}
\label{mobius}
Let  $(\mathcal{S}, \leq)$ be the subgroup lattice of $Aut(G)$ and $\mu(*,*)$ be the M\"obius function of $(\mathcal{S}, \leq)$.  Let $S_0$ be the subgroup of $Aut(G)$ consisting of the identity automorphism of $G$. Then
$$ L(G,k) = N_=(S_0) = \sum_{S \in \mathcal{S}: S \ge S_0} \mu(S_0, S) N_{\ge}(S).$$
\end{theorem}

Using the formula above, the number of $N_\ge(S)$ terms we have to compute to determine $L(G,k)$ is $|\mathcal{S}|$ as opposed to $\Omega( 2^{|Aut(G)|})$  in the PIE formulation. This of course comes at a price -- we must now find all the subgroups of $Aut(G)$, determine the structure of   $(\mathcal{S}, \leq)$, and then compute the $\mu(S_0, S)$ values.  We leave it up to the reader to apply Theorem~\ref{mobius} when $Aut(G) \cong Z_t, D_t, Z_t \times Z_2$ or $D_t \times Z_2$.

Another direction in which we use Lemma~\ref{subgroup} is a generalization of Theorem~\ref{cyclic}.  

\begin{theorem}
\label{cyclic-general}
Let $S^*$ be a subset of $Aut(G)$ that does not contain the identity automorphism.  Suppose every non-trivial subgroup of $G$ contains at least one element of $S^*$.  Then
$$ L(G,k) = \sum_{S \subseteq S^*} (-1)^{|S|} N_{\ge}(S).$$ 
\end{theorem}
 
In the above theorem, the number of $N_\ge(S)$ terms to compute is $2^{|S^*|}$ so the smaller $|S^*|$ is, the better.  Finding the smallest $S^*$, however, is non-trivial; it is the {\it hitting set problem} and is known to be NP-hard in general~\cite{GJ}. Again, we leave it up to the reader to formulate a generalization  of Theorem~\ref{dihedral} similar to Theorem~\ref{cyclic-general}.

Finally, we end with the following open problem.  Let $\mathcal{G}$ be a family of graphs for which there is an efficient algorithm for testing graph isomorphism.  Can the distinguishing number of graphs in $\mathcal{G}$ be computed efficiently?  More specifically, is there a polynomial-time algorithm for computing $D(G)$ if $G$ is a bounded degree graph or a bounded genus graph?

\section*{Acknowledgments}

The second author would like to thank Jeb Willenbring for discussions she had with him on Section 3 of the paper.  

\bibliography{ref}
\bibliographystyle{abbrv}

\end{document}